\documentclass[11pt]{article}
\usepackage{amsfonts,latexsym}

\newcommand\cV{\mathcal{V}}
\newcommand\cL{\mathcal{L}}
\newcommand\cT{\mathcal{T}}
\newcommand\cW{\mathcal{W}}
\newcommand\cU{\mathcal{U}}

\begin{document}
\title{Balanced realizations of discrete-time stable all-pass systems and
 the tangential Schur algorithm}

\author{Bernard Hanzon\thanks{School of Mathematical Sciences, University College, Cork, Ireland, Tel.:
+353-21-4902376, F ax +353-21-4903784, E-mail: b.hanzon@ucc.ie}
   \and Martine Olivi\thanks{Corresponding author. INRIA Sophia-Antipolis,
BP 93, 
06902 Sophia-Antipolis Cedex, FRANCE, 
 Tel.: 33-492-387877, 
Fax: 33 4 92 38 78 58, E-mail: olivi@sophia.inria.fr}  
 \and Ralf L.M. Peeters\thanks{Dept. Mathematics, Universiteit Maastricht,
   P.O. Box 616, 6200 MD 
     Maastricht, The Netherlands, Tel.: 31-43-3883365, Fax: 31-43-3884910,
   E-mail: ralf.peeters@math.unimaas.nl}}
\date{}
\maketitle

\begin{abstract} In this paper, the connections
 are investigated between two different approaches towards the
 parametrization of multivariable stable all-pass systems in
 discrete-time. The first approach involves the tangential Schur
 algorithm, which employs linear fractional transformations. It
 stems from the theory of reproducing kernel Hilbert spaces and
 enables the direct construction of overlapping local
 parametrizations using Schur parameters and interpolation points.
 The second approach proceeds in terms of state-space realizations.
 In the scalar case, a balanced canonical form exists that can also
 be parametrized by Schur parameters. This canonical form can be
 constructed recursively, using unitary matrix operations. Here,
 this procedure is generalized to the multivariable case by
 establishing the connections with the first approach. It gives
 rise to balanced realizations and overlapping canonical forms
 directly in terms of the parameters used in the tangential Schur
 algorithm.
\end{abstract}

\newcommand{\RR}{\mathbb R}
\newcommand{\CC}{\mathbb C}
\newcommand{\UU}{\mathbb U}
\newcommand{\TT}{\mathbb T}
\newcommand{\ZZ}{\mathbb Z}

\newcommand{\eml}{\hfil\break\hbox{}\hfil\break}

\newtheorem{theorem}{Theorem}[section]
\newtheorem{proposition}[theorem]{Proposition}
\newtheorem{lemma}[theorem]{Lemma}
\newtheorem{definition}[theorem]{Definition}
\newtheorem{corollary}[theorem]{Corollary}
\newtheorem{notation}[theorem]{Notation}


\section{Introduction}
\label{intro}%

Stable all-pass systems of finite order have several applications
in linear systems theory. Within the fields of system
identification, approximation and model reduction, they have been
used in connection with the Douglas-Shapiro-Shields factorization,
see e.g., \cite{DSS70,Bar90,L-O98,F-O98}, to obtain effective
algorithms for various purposes, such as $H_2-$model order reduction
of stable linear systems. The class of stable all-pass transfer functions of finite order is
bijectively related to the class of 
rational inner functions, of which the differential structure has
been studied in \cite{ABG94}. There, a parametrization for the
multivariable case has been obtained by means of a recursive
procedure, the {\em tangential Schur algorithm}, which involves
Schur parameter vectors, interpolation points and normalized
direction vectors. In the scalar case, a single coordinate chart
suffices to entirely describe the manifold of stable all-pass (or
lossless) systems of a fixed finite order. In the multivariable
case, the approach leads to atlases of generic charts covering
these manifolds.

In another line of research, balanced state-space canonical forms
have been constructed for various classes of linear systems, with
special properties of these classes (such as stability) built in,
see e.g., \cite{Obe87,H-O98}. Balanced realizations are well-known
to have numerical advantages and are useful for model reduction
purposes in conjunction with balance-and-truncate type procedures.
A balanced canonical form for SISO stable all-pass systems in continuous time was presented
in \cite{OberMunich1987}.
In the constructions of \cite{H-O97}, \cite{H-O98}, the case of stable all-pass
systems in {\em continuous-time} plays a central role. In the
scalar case, the resulting canonical form for lossless systems is balanced with a
positive upper triangular reachability matrix. In the
multivariable case, Kronecker indices and nice selections are used
to arrive at balanced overlapping canonical forms for lossless systems. For {\em
discrete-time} stable all-pass systems, canonical forms can be
obtained from the results in continuous-time by application of a
bilinear transformation. However, this destroys certain nice
properties of the canonical form; e.g., truncation of state
components no longer leads to reduced order systems that are
balanced and in canonical form. Therefore, the ideas of  \cite{H-O97} and 
\cite{H-O98} are applied in \cite{HP2000} to the scalar
discrete-time stable all-pass case directly. This leads to a
balanced canonical form with the desired properties, for which it
turns out that it can in fact be parametrized using Schur
parameters.

In this paper, the connections between these two approaches are
investigated and the results of \cite{HP2000} are generalized to
the multivariable case. Moreover, the atlases developed in
\cite{ABG94} are supplied with balanced state-space realizations
in terms of the parameters used in the tangential Schur algorithm.
In Section \ref{recursive}, we introduce a mapping which acts on
proper rational matrix functions and which has a particularly
simple expression on the level of realizations, see formula
(\ref{statespacerecursion}). When specialized to the case of
lossless functions, it allows for a recursive construction of
balanced realizations that can be implemented as a product of
unitary matrices. This construction has several advantages from a
numerical point of view. In Section \ref{linear}, we recall some
basic results on linear fractional transformations of lossless
functions, which are at the heart of the tangential Schur
algorithm, described in Section \ref{tangential}. The results of
these sections are mostly well-known, except for the nice
factorization of $J$-inner functions in Proposition
\ref{thetafact}, which we could not find in the literature and
which enables an elegant presentation of these results. The main
technical results are presented in Section \ref{connection} and
they describe the complete connection between the construction of
Section \ref{recursive} and the tangential Schur algorithm in its
more general form. In Section \ref{overlapping} the results are
used to construct atlases for the manifold of lossless systems of
a fixed McMillan degree and for the manifold of stable linear
systems of a fixed McMillan degree.

\section{Preliminaries}
\label{preliminaries}%

In this section a brief overview of results from the literature is
presented to provide the necessary background for the main
constructions of this paper. They concern realization theory,
balancing and the basic theory of lossless and $J$-inner
functions.
\\[5mm]
\emph{Realization theory and balancing.}
\\[2mm]
Consider a linear time-invariant state-space system in discrete
time with $m$ inputs and $p$ outputs:
\begin{eqnarray}
  & & x_{t+1} = A x_t + B u_t, \\
  & & y_t     = C x_t + D u_t,
\end{eqnarray}
with $t \in {\mathbb Z}$, $x_t \in {\mathbb C}^n$ for some
nonnegative integer $n$, $u_t \in {\mathbb C}^m$, $y_t \in
{\mathbb C}^p$. Furthermore, the matrices $A$, $B$, $C$ and $D$
with complex-valued entries are of compatible sizes: $n \times n$,
$n \times m$, $p \times n$ and $p\times m$, respectively. The
corresponding transfer matrix of this system is given by
$G(z)=D+C(zI_n-A)^{-1}B$, which is a $p \times m$ matrix with
rational functions as its entries. Two state-space systems are
called input-output equivalent if they have the same transfer
matrix. Here, two rational (matrix) functions are identified if
they agree almost everywhere, i.e., common factors are always
canceled. Note that a transfer matrix in this set-up is always
proper (i.e., analytic at infinity), $D$ being the value at
infinity.
\\[5mm]
From realization theory it follows that, conversely, any $p \times
m$ rational matrix function $G(z)$ analytic at infinity can be
written in the form
\begin{equation}
\label{realization}%
  G(z) = D + C (zI_n-A)^{-1} B
\end{equation}
where $(A,B,C,D)$ is an appropriate quadruple of matrices and $n$
a suitable state space dimension. Such a quadruple with the
associated expression (\ref{realization}) is called a {\em
state-space realization} of $G(z)$. To such a realization we
associate the block-partitioned matrix
\begin{equation}
  R = \left[ \begin{array}{cc} D & C \\ B & A \end{array} \right]
\end{equation}
which we call the {\em realization matrix} as in \cite{HP2000}. It will play an
important role in the sequel.
\\[5mm]
If all possible realizations of $G(z)$ have state-space dimension
at least as big as the state-space dimension $n$ of $(A,B,C,D)$
then the latter realization is said to be minimal, and $n$ is called the
order or the McMillan degree of the transfer function. Two minimal
realizations $(A,B,C,D)$ and
$(A^\prime,B^\prime,C^\prime,D^\prime)$ of a given function $G(z)$
are always similar: there exists a unique invertible matrix $T$
such that
\begin{equation}
\label{similarity}%
  \left[ \begin{array}{cc} D^\prime & C^\prime \\ B^\prime &
  A^\prime \end{array} \right]
  = \left[ \begin{array}{cc} I_p & 0 \\ 0 & T \end{array} \right]
  \left[ \begin{array}{cc} D & C \\ B & A \end{array} \right]
  \left[ \begin{array}{cc} I_m & 0 \\ 0 & T^{-1} \end{array}
  \right].
\end{equation}
As is well-known, an output pair $(C,A)$ is observable if the
observability matrix
\begin{equation}
  \left[ \begin{array}{c} C \\ CA \\ \vdots \\ CA^{n-1}
  \end{array} \right]
\end{equation}
has full column rank. An input pair $(A,B)$ is reachable if the
associated reachability matrix
\begin{equation}
  \left[ \begin{array}{cccc} B & AB & \ldots & A^{n-1}B
  \end{array} \right]
\end{equation}
has full row rank. It is well-known that a realization $(A,B,C,D)$
is minimal if and only if $(C,A)$ is observable and $(A,B)$ is
reachable; cf., e.g., \cite{Kailath}. In this case, the poles of
the function $G(z)$ are the eigenvalues of $A$, and its McMillan
degree is equal to the sum of the degrees (see \cite[Sect.\
6.5]{Kailath}) of all the poles of $G(z)$. This provides an
alternative definition for the McMillan degree, which generalizes
to the larger class of proper and non-proper rational functions by including a possible
pole
at infinity with its appropriate degree. 
\\[5mm]
Let $(A,B,C,D)$ be some realization of a transfer function. If the
eigenvalues of $A$ all belong to the open unit disk, then the
matrix $A$ is called (discrete-time) asymptotically stable, and
$(A,B,C,D)$ an asymptotically stable realization. If $(A,B,C,D)$
is an asymptotically stable realization, then the controllability
Gramian $W_c$ and the observability Gramian $W_o$ are well defined
as the exponentially convergent series
\begin{equation}
\label{Gramians}%
  W_c = \sum_{k=0}^{\infty} A^k B B^* (A^*)^k,
  \hspace{5mm} W_o = \sum_{k=0}^{\infty} (A^*)^k C^* C A^{k},
\end{equation}
where the notation $^*$ is used to denote Hermitian transposition
of a matrix (i.e., the joint action of matrix transposition and
complex conjugation of the matrix entries). The Gramians are
characterized as the unique (and positive semi-definite) solutions
of the respective Lyapunov-Stein equations
\begin{eqnarray}
\label{Steineq1}%
  W_c - A W_c A^* & = & BB^*, \\
\label{Steineq2}%
  W_o - A^* W_o A & = & C^{*}C.
\end{eqnarray}
Moreover, under asymptotic stability of $A$ it holds that $W_c$ is
positive definite if and only if the pair $(A,B)$ is reachable,
and $W_o$ is positive definite if and only if the pair $(C,A)$ is
observable. A minimal and asymptotically stable realization
$(A,B,C,D)$ of a transfer function is called {\em balanced} if its
observability and controllability Gramians are both diagonal and
equal. Any minimal and asymptotically stable realization
$(A,B,C,D)$ is similar to a balanced realization. The concept of
balanced realizations was first introduced in \cite{Moore} in the
continuous time case and used for model reduction. In \cite{PS82}
the same was done for the discrete time case. Balanced
realizations are now a well-established tool which often exhibit
good numerical properties.
\\[5mm]
\emph{$J$-inner, $J$-unitary and lossless functions.}
\\[2mm]
For any matrix function $R(z)$, we define the matrix functions
$R^*(z)$ and $R^{\sharp}(z)$ by
\begin{equation}
\label{Rsharp}%
  R^*(z) := R(\overline{z})^*,
  \hspace{5mm} \mbox{\rm and} \hspace{5mm}
  R^{\sharp}(z) := R^*(z^{-1}).
\end{equation}
Note that if $z$ lies on the unit circle, then $R^{\sharp}(z) =
R(z)^*$. For (square) Hermitian matrices $P$ and $Q$, either of
the notations $P \leq Q$ and $Q \geq P$ will be used to express
that $Q-P$ is a positive semi-definite matrix. As is well-known
this introduces a partial ordering on the set of square Hermitian
matrices.
\\[5mm]
Now let $\Sigma$ be a constant $k \times k$ matrix which is both
Hermitian and unitary. Note that $\Sigma$ is unitarily similar to
a signature matrix, i.e., there exists a unitary $k \times k$
matrix $U$ for which $U \Sigma U^*$ attains the form $\left[
\begin{array}{cc} I_q & 0 \\ 0 & -I_r \end{array}\right]$ for some
non-negative integers $q$ and $r$ with $q+r=k$.
\\[5mm]
A square rational matrix function $\Theta(z)$ of size $k \times k$
is called $\Sigma$-inner (in the unit disk), if at every point of
analyticity $z$ of $\Theta(z)$ it satisfies
\begin{eqnarray}
\label{Jinnerin}%
  \Theta(z)^* \Sigma \Theta(z) & \leq & \Sigma,~~~ |z|<1, \\
\label{Jinneron}%
  \Theta(z)^* \Sigma \Theta(z) & =    & \Sigma,~~~ |z|=1, \\
\label{Jinnerout}
  \Theta(z)^* \Sigma \Theta(z) & \geq & \Sigma,~~~ |z|>1.
\end{eqnarray}
It can be shown that this definition contains redundancy in the
sense that either one of these three defining properties
(\ref{Jinnerin})-(\ref{Jinnerout}) is implied by the other two;
cf., e.g., \cite{Pot60,Dym,Getal}.
\\[5mm]
Any rational matrix function $\Theta(z)$ which satisfies the
property (\ref{Jinneron}) is called $\Sigma$-unitary. For such
functions, replacing $\overline{z}$ by $z^{-1}$ on the unit
circle, the identity (\ref{Jinneron}) extends almost everywhere by
analytic continuation, so that a rational $\Sigma$-unitary
function is invertible and its inverse is given by
\begin{equation}
\label{thetainv}%
  \Theta(z)^{-1} = \Sigma \Theta^{\sharp}(z) \Sigma.
\end{equation}
This applies \emph{a fortiori} to $\Sigma$-inner functions. A
function is called $\Sigma$-lossless if it is $(-\Sigma)$-inner.
If $\Theta(z)$ is $\Sigma$-inner, then $\Theta^*(z)$ is also
$\Sigma$-inner whereas both $\Theta^\sharp(z)$ and
$\Theta(z)^{-1}$ are $\Sigma$-lossless. The class of
$\Sigma$-inner functions is closed under multiplication.
\\[5mm]
As usual, an $I_k$-lossless function is just called {\em
(discrete-time) lossless}, or {\em (discrete-time) stable
all-pass} and an $I_k$-inner function is simply called {\em
inner}. Throughout this paper we will be much concerned with
lossless and inner functions of size $p \times p$, and with
$J$-unitary and $J$-inner functions of size $2p \times 2p$, where
$J$ denotes the following signature matrix:
\begin{equation}
  J = \left[ \begin{array}{cc} I_p & 0 \\ 0 & -I_p \end{array}
  \right].
\end{equation}
Note that $J$-inner and $J$-lossless functions in general may have
poles everywhere in the complex plane, but inner functions are
analytic inside the unit disk and lossless functions are analytic
outside the unit disk (including at the point at infinity).
Therefore a rational lossless function $G(z)$ is proper, having an
inverse which is inner and given by
\begin{equation}
\label{inverselossless}%
  G(z)^{-1} = G^{\sharp}(z).
\end{equation}
A scalar rational function is a {\em Blaschke product} if it maps
the unit circle into the unit circle. It can be shown that such a
function can be written as the product of a finite number of {\em
Blaschke factors}, which are rational functions of first order
that map the unit circle into the unit circle.
\begin{lemma}
\label{degreelossless}%
Let $G(z)$ be a $p \times p$ lossless function of McMillan degree
$n$. Then its determinant $\det G(z)$ is a Blaschke product. If
one writes it in the form $\det G(z) = e(z)/d(z)$ with $e$ and $d$
co-prime polynomials, then the degree of $d(z)$ is equal to $n$.
\end{lemma}
{\em Proof.~}%
See e.g. \cite{BOiuf}.
\hfill $\Box$%
\\[5mm]

 We finally describe some particular $J$-inner functions
that will be intensively used in the sequel, namely the constant $J$-unitary
matrices  and the $J$-inner functions of McMillan degree one. Details and
proofs can be found in \cite{Dym}.

Every \emph{constant} $J$-unitary matrix can be represented in a
unique way (see \cite[Thm.\ 1.2]{Dym}) as follows:
\begin{equation}
\label{Junitaryrep}%
  M = H(E) \left[ \begin{array}{cc} P & 0 \\ 0 & Q \end{array}
  \right],
\end{equation}
where $P$ and $Q$ are $p \times p$ unitary matrices and $H(E)$
denotes the Halmos extension of a strictly contractive $p \times
p$ matrix $E$ (i.e., such that $I-E^*E>0$). This Halmos extension
$H(E)$ is defined by
\begin{equation}
\label{Halmos}
  H(E) = \left[ \begin{array}{cc}
    (I-EE^*)^{-1/2}     & E (I-E^*E)^{-1/2} \\
    E^* (I-EE^*)^{-1/2} & (I-E^*E)^{-1/2}
  \end{array} \right] = \left[ \begin{array}{cc}
    (I-EE^*)^{-1/2}     & (I-EE^*)^{-1/2} E \\
    (I-E^*E)^{-1/2}E^*  & (I-E^*E)^{-1/2}
  \end{array} \right].
\end{equation}
It holds that $H(E)$ is Hermitian, $J$-unitary and invertible with
inverse $H(E)^{-1}=H(-E)$. Also, $K H(E) K = H(E^*)$ for $K =
\left[ \begin{array}{cc} 0 & I_p \\ I_p & 0 \end{array} \right]$.
\\[5mm]
Two other important forms of $J$-unitary matrices are described in the
following lemma ( see \cite[Thm.1.3]{Dym}). 

\begin{lemma}
\label{XxYx}
If $x$ is a
$2p$-vector such that 
$x^* J x \neq 0$, then the set of $2p \times 2p$ matrices of the
form
\begin{equation}
\label{xcheck}%
  X_x(\alpha) = I_{2p} + (\alpha-1) x (x^*Jx)^{-1} x^* J,
  \hspace{5mm} \alpha \in {\mathbb{C}},
\end{equation}
is closed under multiplication. We have that
\begin{equation}
\label{xcheckprop}%
  X_x(\alpha) X_x(\beta) = X_x(\alpha\beta),
  \hspace{5mm} \det X_x(\alpha) = \alpha,
  \hspace{5mm} X_x(\alpha)^{-1} = X_x(\alpha^{-1}),
  \hspace{2mm} \mbox{\rm for~} \alpha \neq 0,
\end{equation}
\begin{equation}
J- X_x(\alpha)J X_x(\beta)^*=(1-\alpha\beta^*) x(x^*Jx)^{-1}x^*J,
\end{equation}
and $X_x(\alpha)$ is $J$-unitary if and only if $|\alpha|=1$.
\\[5mm]
If $x$ is a $2p$-vector such that $x^* J x = 0$, then the
set of matrices of the form
\begin{equation}
\label{xdot}%
  Y_x(\alpha) = I_{2p} + \alpha x x^* J,
  \hspace{5mm} \alpha \in {\mathbb{C}},
\end{equation}
is closed under multiplication. We have that
\begin{equation}
\label{xdotprop}%
  Y_x(\alpha) Y_x(\beta) = Y_x(\alpha+\beta),
  \hspace{5mm} \det Y_x(\alpha) = 1,
  \hspace{5mm} Y_x(\alpha)^{-1} = Y_x(-\alpha),
\end{equation}
\begin{equation}
J- Y_x(\alpha)J Y_x(\beta)^*=(\alpha+\beta^*) x x^*,
\end{equation}
and $Y_x(\alpha)$  
is $J$-unitary if and only if $\alpha$ is pure imaginary.
\end{lemma}

We now come to the description of 
$J$-inner matrix functions of McMillan degree $1$.  
 We shall call
these functions `elementary $J$-inner functions' or `elementary
$J$-inner factors', since it has been shown by Potapov (see
\cite{Pot60}) that every $J$-inner matrix function of McMillan
degree $n$ can be decomposed into a product of $n$ such elementary
$J$-inner factors. Elementary $J$-inner functions are obtained by 
applying the functions $X_x$ and $Y_x$ to appropriate
scalar functions of degree one, namely the Blaschke factor $b_w(z)$ defined by
\begin{equation}
  b_w(z) := \frac{z-w}{1-\overline{w}z},~~~ w\notin \TT
\end{equation}
and the Carath\'eodory function,
\begin{equation}
  c_w(z) := \frac{z+w}{z-w},~~~ w\in \TT,
\end{equation}
where $\TT$ denotes the unit circle. 
\begin{theorem}
\label{Jinnerfactors}
 Let $\Phi(z)$ be an elementary $J$-inner factor with a pole at
 $z=1/\overline{w}$. Then,  apart from  a constant $J$-unitary multiplier on 
 the right, $\Phi(z)$ must be in one of the following three forms:
\item[(1)]If $w\in \CC$, $w \notin \TT$,
\begin{equation}
\label{XBlaschke}%
  \Phi(z) = X_x(b_w(z)) = I_{2p} +
  (b_w(z)-1) \frac{x x^* J}{x^* J 
  x},
\end{equation}
(see Eqn.\ (\ref{xcheck})) for some  $x \in {\mathbb{C}}^{2p }$ such that : 
\begin{equation}
\label{P>0}
\frac{x^* J x}{1-|w|^2}>0.
\end{equation}
\item[(2)] If $w\in\TT$,
\begin{equation}
\label{Ydegone}%
  \Phi(z) = Y_x(-\delta\ c_w(z))
          = I_{2p}-\delta\ c_w(z) \ x x^* J,
 \end{equation}
(see Eqn.\ (\ref{xdot})) for some nonzero $x \in \CC^{2p}$ 
such that $x^*Jx=0$, and for some real $\delta>0$.
\item[(3)] If $w=\infty$, then 
\begin{equation}
\label{XBlaschkeinfty}%
  \Phi(z) = X_x\left(\frac{1}{z}\right) = I_{2p} + \left(\frac{1}{z}-1\right) \frac{x x^* J}{x^* J
  x},
\end{equation}
for some  $x \in {\mathbb{C}}^{2p }$ such that $x^* J x<0$.
\end{theorem}


\section{Recursive construction of balanced realizations of lossless systems}
\label{recursive}%
In this section we will present a number of results that will be used
in the proposed recursive  
construction of balanced state-space realizations of lossless
functions.
\\[1mm]
   With each pair $(U,V)$ of $(p+1) \times (p+1)$ matrices we
associate a mapping ${\cal F}_{U,V}$ which is defined to act on
proper rational $p \times p$ matrix functions $G(z)$ as follows.
\begin{equation}
\label{FUV}%
  {\cal F}_{U,V}:
  G(z) \mapsto F_1(z) + \frac{F_2(z)F_3(z)}{z-F_4(z)},
\end{equation}
with $F_1(z)$ of size $p \times p$, $F_2(z)$ of size $p \times 1$,
$F_3(z)$ of size $1 \times p$ and $F_4(z)$ scalar, defined by:
\begin{equation}
\label{matF}%
  F(z) = \left[ \begin{array}{cc} F_1(z) & F_2(z) \\
  F_3(z) & F_4(z) \end{array} \right] =
  V \left[ \begin{array}{cc} 1 & 0 \\ 0 & G(z) \end{array}
  \right] U^*.
\end{equation}
A state-space realization of $\widetilde{G}(z) = {\cal
F}_{U,V}(G(z))$ can be obtained by working directly on any
realization matrix $R$ of $G(z)$. This is the content of the
following proposition.
\begin{proposition}
\label{ssrecursion}%
Let $G(z)$ be a $p \times p$ proper rational transfer function and
$(U,V)$ a pair of $(p+1)\times (p+1)$ matrices, then
$\widetilde{G}(z) = {\cal F}_{U,V}(G(z))$ is well-defined. Let
$(A,B,C,D)$ be a state-space realization of $G(z)$ with
$n$-dimensional state-space. Then a state-space realization
$(\widetilde{A},\widetilde{B},\widetilde{C},\widetilde{D})$ of
$\widetilde{G}(z)$ with $(n+1)$-dimensional state-space is given
by:
\begin{equation}
\label{statespacerecursion}%
  \left[ \begin{array}{cc} \widetilde{D} & \widetilde{C} \\
  \widetilde{B} & \widetilde{A} \end{array} \right] =
  \left[ \begin{array}{cc} V & 0 \\ 0 & I_n \end{array} \right]
  \left[ \begin{array}{ccc} 1 & 0 & 0 \\ 0 & D & C \\ 0 & B & A
  \end{array} \right] \left[ \begin{array}{cc} U^* & 0 \\ 0 & I_n
  \end{array} \right].
\end{equation}
\end{proposition}
{\em Proof.~}%
Since $G(z)$ is proper, $F_4(z)$ is proper too and therefore
$z-F_4(z)$ does not vanish identically so that ${\cal
F}_{U,V}(G(z))$ is well-defined.
\\[5mm]
Now observe that the right-hand side of Eqn.\
(\ref{statespacerecursion}) provides a realization matrix with
$n$-dimensional state-space for the transfer function $F(z)$ given
by Eqn.\ (\ref{matF}):
\[
  \left[ \begin{array}{cc} D_F & C_F \\ B_F & A_F \end{array}
  \right] = \left[ \begin{array}{cc} V & 0 \\ 0 & I_n \end{array}
  \right] \left[ \begin{array}{ccc} 1 & 0 & 0 \\ 0 & D & C \\
  0 & B & A \end{array} \right] \left[ \begin{array}{cc} U^* & 0
  \\ 0 & I_n \end{array} \right].
\]
Here $D_F$ is $(p+1) \times (p+1)$, $C_F$ is $(p+1) \times n$,
$B_F$ is $n \times(p+1)$ and $A_F=A$ is $n \times n$. We shall
prove that changing the partitioning of this matrix to the one in
(\ref{statespacerecursion})
\[
  \left[ \begin{array}{cc} D_F & C_F \\ B_F & A_F \end{array}
  \right] = \left[ \begin{array}{cc} \widetilde{D} &
  \widetilde{C} \\ \widetilde{B} & \widetilde{A} \end{array}
  \right],
\]
where $\widetilde{D}$ is $p \times p$, $\widetilde{C}$ is $p
\times (n+1)$, $\widetilde{B}$ is $(n+1) \times p$ and
$\widetilde{A}$ is $(n+1) \times (n+1)$, gives a realization of
$\widetilde{G}(z)$ with $n+1$-dimensional state-space. To see
this, let the matrices $\widetilde{C}$, $\widetilde{B}$ and
$\widetilde{A}$ be partitioned as follows:
\[
  \left[ \begin{array}{cc} \widetilde{D} & \widetilde{C} \\
  \widetilde{B} & \widetilde{A} \end{array} \right]
  = \left[ \begin{array}{c|cc} \widetilde{D} &
  \widetilde{C}_1 & \widetilde{C}_2 \\
  \hline
  \widetilde{B}_1 &\widetilde{A}_1 & \widetilde{A}_2 \\
  \widetilde{B}_2 & \widetilde{A}_3 & A_F
  \end{array}\right]
  = \left[ \begin{array}{cc|c} \widetilde{D} &
  \widetilde{C}_1 & \widetilde{C}_2 \\
  \widetilde{B}_1 &\widetilde{A}_1 & \widetilde{A}_2 \\
  \hline
  \widetilde{B}_2 & \widetilde{A}_3 & A_F \end{array} \right]
  = \left[ \begin{array}{cc} D_F & C_F \\
  B_F & A_F \end{array} \right],
\]
where $\widetilde{C}_1$ is $p \times 1$, $\widetilde{C}_2$ is $p
\times n$, $\widetilde{B}_1$ is $1 \times p$, $\widetilde{B}_2$ is
$n \times p$, $\widetilde{A}_2$ is $1 \times n$, $\widetilde{A}_3$
is $n \times 1$ and $\widetilde{A}_1$ is scalar. Note that it
follows that $F_4(z)$ has the realization
$(A_F,\widetilde{A}_3,\widetilde{A}_2,\widetilde{A}_1)$ and
therefore
\[
  z - \widetilde{A}_1 - \widetilde{A}_2 (z I_{n}-A_F)^{-1}
  \widetilde{A}_3 = z-F_4(z).
\]
Using this together with the well-known formula for the inverse of
a partitioned matrix (see \cite{Dym} or \cite{Kailath}) we compute
\[
  (z I_{n+1} -\widetilde{A})^{-1} =
  \left[ \begin{array}{cc} (z-F_4(z))^{-1} & - (z-F_4(z))^{-1}
  \widetilde{A}_2 (z I_n -A_F)^{-1} \\
  -(z I_n-A_F)^{-1} \widetilde{A}_3 (z-F_4(z))^{-1} & \Delta(z)
  \end{array} \right]
\]
with
\[
  \Delta(z) = (z I_n -A_F)^{-1} + (z I_n-A_F)^{-1}
  \widetilde{A}_3 (z-F_4(z))^{-1} \widetilde{A}_2
  (z I_n-A_F)^{-1}.
\]
We then have
\begin{eqnarray*}
  \widetilde{D} + \widetilde{C} (z I_{n+1} - \widetilde{A})^{-1}
  \widetilde{B} & = & \widetilde{D} + \left[ \begin{array}{cc}
  \widetilde{C}_1 & \widetilde{C}_2 \end{array} \right]
  (z I_{n+1} - \widetilde{A})^{-1} \left[ \begin{array}{c}
  \widetilde{B}_1 \\ \widetilde{B}_2 \end{array} \right] \\
  & = & \widetilde{D} + \widetilde{C}_1 (z-F_4(z))^{-1}
  \widetilde{B}_1 + \widetilde{C}_1  (z-F_4(z))^{-1}
  \widetilde{A}_2 (z I_n-A_F)^{-1} \widetilde{B}_2 + \\
  & & + \widetilde{C}_2 (z I_n-A_F)^{-1} \widetilde{A}_3
  (z-F_4(z))^{-1} \widetilde{B}_1 + \widetilde{C}_2
  (z I_n-A_F)^{-1} \widetilde{B}_2 + \\
  & & + \widetilde{C}_2 (z I_n-A_F)^{-1} \widetilde{A}_3
  (z-F_4(z))^{-1} \widetilde{A}_2 (z I_n-A_F)^{-1} \widetilde{B}_2
  \\
  & = & F_1(z) + F_2(z) (z-F_4(z))^{-1} F_3(z),
\end{eqnarray*}
which finally proves that
$(\widetilde{A},\widetilde{B},\widetilde{C},\widetilde{D})$ is a
state-space realization of $\widetilde{G}(z)$ with
$(n+1)$-dimensional state-space.
\hfill $\Box$%
\\[5mm]
The following proposition characterizes the nature of minimal
balanced state-space realizations, in discrete time, of rational
lossless functions. The fact that a lossless function admits a
minimal balanced realization for which the realization matrix is
unitary is well-known from the literature; see, e.g.,
\cite{Getal}. However, the converse result, which asserts that
unitary realization matrices correspond to possibly non-minimal
realizations of lossless functions, appears to be novel and
generalizes \cite[Prop.\ 2.4]{HP2000} to the multivariable case.
\begin{proposition}
\label{unitaryrealization}%
(i) For any minimal balanced realization of a $p \times p$
rational lossless function the observability and controllability
Gramians are both equal to the identity matrix and the associated
realization matrix is unitary.
\\
(ii) Conversely, if the realization matrix associated with a
realization $(A,B,C,D)$ of order $n$ of some $p \times p$ rational
function $G$ is unitary, then $G$ is lossless of McMillan degree
$\leq n$. The realization is minimal if and only if $A$ is
asymptotically stable and then it is balanced.
\end{proposition}
{\em Proof.~}%
(i) As is well known, any minimal realization $(A,B,C,D)$ of a
lossless function $G$ is asymptotically stable and satisfies (see
\cite{Getal})
\begin{equation}
\label{real-lossless-eq}%
  \left[ \begin{array}{cc} D & C \\ B & A \end{array} \right]
  \left[ \begin{array}{cc} I_p & 0 \\ 0 & P \end{array} \right]
  \left[ \begin{array}{cc} D^* & B^* \\ C^* & A^* \end{array} \right]
  = \left[ \begin{array}{cc} I_p & 0 \\ 0 & P \end{array} \right],
\end{equation}
for some unique positive definite $n \times n$ matrix $P$ which is
precisely the controllability Gramian. Since $P$ is positive
definite, it can be reduced to identity by using any state space
transformation matrix $T$ such that $P=TT^*$, producing a new
minimal state space realization
$(\widetilde{A},\widetilde{B},\widetilde{C},\widetilde{D}) =
(T^{-1}AT,T^{-1}B,CT,D)$. Then Eqn.\ (\ref{real-lossless-eq})
asserts that the associated new realization matrix is unitary.
Upon changing the order of multiplication of this unitary
realization matrix and its Hermitian transpose it follows that the
new observability Gramian is also equal to identity. Hence,
$(\widetilde{A},\widetilde{B},\widetilde{C},\widetilde{D})$ is
balanced. Any other minimal balanced realization of $G(z)$ can be
obtained from
$(\widetilde{A},\widetilde{B},\widetilde{C},\widetilde{D})$ by
means of a unitary change of basis of the state space, which
leaves the associated realization matrix unitary.
\\
(ii) Conversely, if a realization $(A,B,C,D)$ has unitary realization matrix,
a standard straightforward calculation shows that
\begin{equation}
  I_p - G(z) G(w)^*  =  (z\overline{w}-1) C (zI_n-A)^{-1}
  (\overline{w} I_n-A^*)^{-1} C^*,
\end{equation}
for every choice of $z$ and $\bar w$ outside the spectrum of $A$.
Setting $w=z$, it follows that $G$ is lossless. Therefore, $G$ has
all its poles inside the open unit disk, which makes asymptotic
stability of $A$ a necessary condition for minimality of
$(A,B,C,D)$. On the other hand asymptotic stability of $A$ implies
that minimality of $(A,B,C,D)$ is equivalent to the existence of
positive definite solutions to the Lyapunov-Stein equations
(\ref{Steineq1})-(\ref{Steineq2}). Unitarity of the realization
matrix, however, implies that $I_n$ is always a positive definite
solution to both these equations, showing asymptotic stability of
$A$ to be a sufficient condition for minimality of $(A,B,C,D)$ as
well. Clearly in such a case $(A,B,C,D)$ is also balanced.
\hfill $\Box$%
\\[5mm]%
{\em Remark.~}%
 Let $G_1$ and $G_2$ be two lossless
functions of McMillan degree $n_1$ and $n_2$, respectively, with
minimal realizations $(A_1,B_1,C_1,D_1)$ and $(A_2,B_2,C_2,D_2)$
both having unitary associated realization matrices. Then the
cascade realization $(A,B,C,D)$ of the lossless function $G=G_1
G_2$ is obtained as
\[
  \left[ \begin{array}{c|c} D & C \\ \hline
    B & A \end{array} \right]
  = \left[ \begin{array}{c|cc} D_1 & 0 & C_1 \\ \hline
    0 & I_{n_2} & 0 \\ B_1 & 0 & A_1 \end{array} \right]
  \left[ \begin{array}{c|cc} D_2 & C_2 & 0 \\ \hline
    B_2 & A_2 & 0 \\ 0 & 0 & I_{n_1} \end{array} \right]
  = \left[ \begin{array}{c|cc} D_1 D_2 & D_1 C_2 & C_1 \\ \hline
    B_2 & A_2 & 0 \\ B_1 D_2 & B_1 C_2 & A_1 \end{array} \right]
\]
which yields again a unitary realization matrix which is minimal.

Three properties of the class of mappings ${\cal F}_{U,V}$ are
collected in the following lemma for later use. These properties
follow straightforwardly from the definition of the mapping ${\cal
F}_{U,V}$ with $U$, $V$ arbitrary $(p+1) \times (p+1)$ matrices.
\begin{lemma}
\label{3properties}%
Let $U$ and $V$ be arbitrary $(p+1) \times (p+1)$ matrices and $P$
and $Q$ both $p \times p$ matrices. Then we have:
\vspace{-0.5cm}
\begin{description}
\item[(i)] ${\cal F}_{U \left[ \begin{array}{cc} 1 & 0 \\
0 & P \end{array} \right], V \left[ \begin{array}{cc} 1 & 0 \\
0 & Q \end{array} \right]}(G(z)) = {\cal F}_{U,V}(QG(z)P^*)$,
\item[(ii)] $\widetilde{G}(z) = {\cal F}_{U,V}(G(z))
\Leftrightarrow \widetilde{G}^*(z) = {\cal F}_{V,U}(G^*(z))$,
\item[(iii)] ${\cal F}_{\left[ \begin{array}{cc} I_p & 0 \\
0 & \xi \end{array} \right] U, \left[ \begin{array}{cc} I_p & 0 \\
0 & \xi \end{array} \right] V}(G(z)) = {\cal F}_{U,V}(G(z))$, for
all $\xi \in {\mathbb C}$ with $|\xi|=1$.
\end{description}
\end{lemma}
The next property shows some of the importance of the ${\cal
F}_{U,V}$ mappings in relation to lossless systems: if $U$ and $V$
are unitary matrices then lossless systems are mapped into
lossless systems. The McMillan degree increases by at most one
under this mapping.
\begin{proposition}
\label{FUVonlossless}%
Let $U$ and $V$ be unitary $(p+1) \times (p+1)$ matrices. \\
The mapping ${\cal F}_{U,V}$ sends a lossless function of order
$n$ to a lossless function of order $\leq n+1$.
\end{proposition}
{\em Proof.~}%
Since $G(z)$ is lossless of order $n$ it has a minimal realization
with a unitary realization matrix of size $n+p$. Formula
(\ref{statespacerecursion}) gives a realization of
$\widetilde{G}={\cal F}_{U,V}(G)$ which is unitary of size
$n+p+1$. Therefore, $\widetilde G(z)$ is a lossless function of
order at most $n+1$.
\hfill $\Box$%
\\[5mm]
Starting with a $p \times p$ unitary matrix $G_0$, interpreted as
a lossless transfer matrix of order zero, and applying the
previous proposition repeatedly, one obtains a class of lossless
transfer matrices and corresponding realizations. The class that
is obtained after $n$ such recursion steps is parametrized by the
finite sequence of unitary matrices $U_1, V_1, U_2, V_2, \ldots,
U_n, V_n$, together with the unitary matrix $G_0$. It consists of
lossless systems of order $\leq n+1$. In Section \ref{connection}
the question which lossless transfer matrix functions can be
obtained in this way will be treated. This generalizes the results
of \cite{HP2000} for the SISO case, where $U_k=I_2$ and $V_k$ is a
$2 \times 2$ Householder reflection matrix, for $k=1,\ldots,n$.


\section{Linear fractional transformations}
\label{linear}%

Linear fractional transformations occur extensively in
representation formulas for the solution of various interpolation
problems \cite{BGR}. They are at the heart of the parametrization
of lossless functions through the tangential Schur algorithm. The
properties of a linear fractional transformation of matrices were
studied in a form adapted to the needs of $J$-theory in
\cite{Pot88} and in an algebraic setting in \cite{Young}. For more
details on and certain proofs of the basic results presented in this
section, we refer to these papers.

Let $\Theta \in {\mathbb{C}}^{2p \times 2p}(z)$ be an invertible
rational matrix in the variable $z$, block partitioned as $\Theta
= \left[ \begin{array}{cc} \Theta_1 & \Theta_2 \\ \Theta_3 &
\Theta_4 \end{array} \right]$ with blocks $\Theta_i$
($i=1,\ldots,4$) of size $p \times p$.

Associated with $\Theta$, let the linear fractional
transformations ${\cal T}_{\Theta}$ and $\widehat{\cal
T}_{\Theta}$ be defined to act on rational matrices $G \in
{\mathbb C}^{p \times p}(z)$ in the following way:
\begin{eqnarray}
  & & {\cal T}_{\Theta}: G \mapsto
  (\Theta_4 G + \Theta_3) (\Theta_2 G + \Theta_1)^{-1}, \\
  & & \widehat{\cal T}_{\Theta}: G \mapsto
  (G \Theta_2 + \Theta_4)^{-1} (G \Theta_1 + \Theta_3).
\end{eqnarray}
The domains of these mappings are denoted by ${\cal M}_{\Theta}$
and $\widehat{\cal M}_{\Theta}$, respectively. They consist of
those $G$ for which the expressions $\Theta_2 G + \Theta_1$ and $G
\Theta_2 + \Theta_4$, respectively, have full rank $p$.  Due to the
invertibility assumption on $\Theta$, it follows that these domains
are residual subsets (countable intersections of open dense subsets) of
${\mathbb   C}^{p \times p}(z)$.  Note that a
linear fractional transformation is fully determined if it is 
specified on such a residual set. 

It is easily established that the following group properties hold
for the LFTs associated with two invertible matrices $\Phi$ and
$\Psi$:
\begin{eqnarray}
\label{ttheta}%
  & & {\cal T}_{\Phi} \circ {\cal T}_{\Psi} =
      {\cal T}_{\Phi\Psi}, \\
\label{tthetahat}%
  & & \widehat{\cal T}_{\Phi} \circ \widehat{\cal T}_{\Psi} =
      \widehat{\cal T}_{\Psi\Phi}.
\end{eqnarray}
Here it may be noted that the compositions of mappings on the left
hand sides of these equations may formally have domains which are
strict subsets of the domains of the mappings on the right hand
sides. However, since the domains of the compositions of mappings
on the left-hand sides are residual sets too, this issue is resolved by
identifying (compositions of) LFTs if they agree on the
intersection of their domains, which will be our policy throughout
this paper.

It then holds that the LFTs ${\cal T}_{\Theta}$ and $\widehat{\cal
T}_{\Theta}$ are bijections from their domains ${\cal M}_{\Theta}$
and $\widehat{\cal M}_{\Theta}$ to their co-domains ${\cal
M}_{\Theta^{-1}}$ and $\widehat{\cal M}_{\Theta^{-1}}$, and that
every LFT can be represented in each of the two forms described
above. We have the following lemma.
\begin{lemma}
\label{That}%
Let $\Theta \in {\mathbb{C}}^{2p \times 2p}(z)$ be an invertible
rational matrix. Then:
\\
(i) ${\cal T}_{\Theta}$ and $\widehat{\cal T}_{\Theta}$ are
bijections. Their inverses are given by
\begin{eqnarray*}
  & & {\cal T}_{\Theta}^{-1} = {\cal T}_{\Theta^{-1}}=
  \widehat{\cal T}_{J \Theta J}, \\
  & & \widehat{\cal T}_{\Theta}^{-1} =
  \widehat{\cal T}_{\Theta^{-1}}= {\cal T}_{J \Theta J}.
\end{eqnarray*}
(ii) ${\cal T}_{\Theta} = \widehat{\cal T}_{J \Theta^{-1} J}$ and
$\widehat{\cal T}_{\Theta} = {\cal T}_{J \Theta^{-1} J}$.
\end{lemma}
We now proceed to study the conditions under which two LFTs
coincide. From the literature, the following result is well known
\cite{Young}.
\begin{lemma}
\label{LFTequiv}%
Let $\Phi, \Psi  \in {\mathbb{C}}^{2p \times 2p}(z)$ be two
invertible rational matrices. It holds that ${\cal T}_{\Phi} =
{\cal T}_{\Psi}$ if and only if there exists a scalar function
$\lambda \in {\mathbb C}(z)$ for which $\Phi = \lambda \Psi$.
\end{lemma}
For our purposes, however, it will be convenient to develop a
slightly specialized version of this lemma, which states that the
same result remains to hold true if two LFTs ${\cal T}_{\Phi}$ and
${\cal T}_{\Psi}$ are merely required to coincide on the subset of
$p \times p$ lossless functions.
\begin{lemma}
\label{LFTlosslessequiv}%
Let $\Phi, \Psi \in {\mathbb{C}}^{2p \times 2p}(z)$ be two
invertible rational matrices. If the $p \times p$ lossless
functions are all contained in ${\cal M}_{\Phi} \cap {\cal
M}_{\Psi}$, and if for all lossless $G$ it holds that ${\cal
T}_{\Phi}(G) = {\cal T}_{\Psi}(G)$, then there exists a scalar
function $\lambda \in {\mathbb C}(z)$ for which $\Phi = \lambda
\Psi$.
\end{lemma}
{\em Proof.~}%
Note that upon application of the bijection ${\cal
T}_{\Psi^{-1}}$, we have that the identity ${\cal
T}_{\Phi}(G)={\cal T}_{\Psi}(G)$ becomes equivalent to ${\cal
T}_{\Psi^{-1}\Phi}(G)=G$. Denoting $\Theta = \Psi^{-1}\Phi$, this
can be rewritten equivalently as $\Theta_4 G+\Theta_3 = G \Theta_2
G + G\Theta_1$, since $(\Theta_2 G + \Theta_1)^{-1}$ exists (for
all lossless $G$). If $G$ is lossless, then also $-G$ is lossless.
Comparing the two expressions which follow from substitution of
these values for $G$, it is found that for all lossless $G$
\begin{eqnarray*}
  & & \Theta_4 G = G \Theta_1, \\
  & & \Theta_3 = G \Theta_2 G.
\end{eqnarray*}
Substitution of $G=I_p$ yields $\Theta_1=\Theta_4$ and
$\Theta_2=\Theta_3$. Setting $G=z^{-1}I_p$, which is again
lossless, yields $\Theta_3=z^{-2}\Theta_2$. In combination with
$\Theta_2=\Theta_3$ it follows that $\Theta_2=\Theta_3=0$.

We are left with the identity $\Theta_1 G = G \Theta_1$ for all
lossless $G$. Note that every unitary matrix is lossless, so that
$G$ can be varied over all the products of signature matrices and
permutation matrices. From this it easily follows that $\Theta_1$
is a scalar multiple of the identity matrix. So there exists a
scalar function $\lambda$ such that $\Theta_1 = \Theta_4 = \lambda
I_p$. (Rationality of $\lambda$ is obvious from rationality of
$\Theta$, which is a consequence of rationality of $\Phi$ and
$\Psi$.) Hence $\Theta = \lambda I_{2p}$, and therefore
$\Phi=\lambda \Psi$.
\hfill $\Box$%
\\[5mm]
{\em Remark.~}%
The proof of Lemma \ref{LFTlosslessequiv} is such that it also
applies to a setting in which all the rational matrix functions
are restricted to be {\em real}. In the {\em complex} case, one
could replace the choice $G=z^{-1}I_p$ (which is lossless of
McMillan degree $p$) by the choice $G=iI_p$, so that the whole
proof involves values for $G$ from the set of unitary matrices
only (i.e., $G$ lossless of McMillan degree $0$). In the real case
with $p>1$ a more careful analysis shows that a proof can also be
designed using orthogonal matrices only. However, in the real case
with $p=1$, the only lossless functions of McMillan degree $0$ are
$G=\pm 1$. Here, lossless functions of McMillan degree $\geq 1$
are needed to obtain the desired result.
\\[5mm]
For our purposes, it is important to study in some more detail the
situation where the $2p \times 2p$ rational matrix $\Theta$ is
$J$-inner and the $p \times p$ rational matrices $G$ on which the
LFT ${\cal T}_{\Theta}$ acts are lossless.
\begin{proposition}
\label{LFTonlossless}%
Let $\Theta$ be $J$-inner of size $2p \times 2p$ and of McMillan
degree $m$. If $G$ is $p \times p$ lossless and of McMillan degree
$n$, then $G \in {\cal M}_{\Theta}$ and the matrix function
$\widehat{G}={\cal T}_{\Theta}(G)$ is also lossless and of
McMillan degree $\leq n+m$.
\end{proposition}
{\em Proof.~}%
 The proof can be easily adapted from \cite[sect.18.2.]{BGR}. 
A realization of $\widehat{G}$ with $(n+m)$-dimensional state-space can be
easily computed (see formula (4.85) in \cite{Kimura}), which proves that  the McMillan degree of $\widehat{G}$ is at most
$n+m$.
\hfill $\Box$%
\\[5mm]
Also, if the LFT corresponds to a $J$-inner matrix, we have the
following elementary results.
\begin{lemma}
\label{LFTproperties}%
Let $\Theta(z)$ be $J$-inner of size $2p \times 2p$. \\
(a) The mappings ${\cal T}_{\Theta}$ and $\widehat{\cal
T}_{\Theta^{\sharp}}$ coincide. \\
(b) If $P$ and $Q$ are $p \times p$ unitary matrices, then
$\Theta_{P,Q}(z) := \Theta(z) \left[ \begin{array}{cc} P & 0 \\
0 & Q \end{array} \right]$ is again $J$-inner, and for all $G \in
{\cal M}_{\Theta_{P,Q}}$:
\[
 {\cal T}_{\Theta_{P,Q}(z)}(G(z))={\cal T}_{\Theta(z)}(QG(z)P^*).
\]
(c) Let $K = \left[ \begin{array}{cc} 0 & I_p \\ I_p & 0
\end{array} \right]$. Then the function defined by
\begin{equation}
\label{analyticoutside}%
  \Theta^o(z) := K \Theta(z^{-1}) K
\end{equation}
is again $J$-inner, and for all $G \in {\cal M}_{\Theta}$ it holds
that ${\cal T}_\Theta(G)^* = {\cal T}_{\Theta^o}(G^*)$.
\end{lemma}
{\em Proof.~}%
(a) Since $\Theta$ is $J$-inner, it holds that $\widehat{\cal
T}_{\Theta^{\sharp}} = \widehat{\cal T}_{J \Theta^{-1}J} = {\cal
T}_{\Theta}$.
\\
(b) The proof that $\Theta_{P,Q}(z)$ is again $J$-inner is
elementary. The rest of this statement follows from the group
property of LFTs and from the simple action of an LFT associated
with a block-diagonal matrix.
\\
(c) Note that ${\cal T}_{\Theta}(G)^* =
(G^*\Theta_2^*+\Theta_1^*)^{-1} (G^*\Theta_4^*+\Theta_3^*) =
\widehat{\cal T}_{K \Theta^* K}(G^*) = {\cal
T}_{JK(\Theta^*)^{-1}KJ}(G^*) = {\cal
T}_{K(\Theta^*)^{\sharp}K}(G^*) = {\cal T}_{\Theta^o}(G^*)$. The
proof that $\Theta^o(z)$ is again $J$-inner is elementary upon
observing that $KJ=-JK$.
\hfill $\Box$%
\\[5mm]
A particular situation of interest occurs when the $J$-inner
matrix function $\Theta(z)$ happens to be of McMillan degree $0$.
In this case $\Theta(z)=M$ is a constant $J$-unitary matrix and
the associated class of LFTs is the one of {\em generalized
M\"obius transformations} ${\cal T}_M$. The following theorem
indicates how the action of such transformations ${\cal T}_M$ on
{\em lossless} functions $G(z)$ can be represented entirely in
terms of {\em balanced} state-space realizations.
\begin{theorem}
\label{moebius}%
The linear fractional transformation ${\cal T}_M$ associated with
a constant $J$-unitary matrix $M$ is a bijection on the set of
lossless functions which preserves the McMillan degree.
\\
Let $M = \left[ \begin{array}{cc} M_1 & M_2 \\ M_3 & M_4
\end{array} \right]$ be the block decomposition of $M$ and let
$G(z)$ be a lossless function with minimal state-space realization
$(A,B,C,D)$. Then, a minimal state-space realization
$(\widetilde{A},\widetilde{B},\widetilde{C},\widetilde{D})$ for
$\widetilde{G}(z)= {\cal T}_M(G(z))$ is given by:
\begin{equation}
\label{realmoebius}%
  \left\{\begin{array}{rcl}
    \widetilde{A}& =& A - B (M_2 D + M_1)^{-1} M_2 C, \\
    \widetilde{B}& =& B (M_2 D + M_1)^{-1}, \\
    \widetilde{C}& =& [M_4 - (M_4 D + M_3)
                       (M_2 D + M_1)^{-1} M_2] C, \\
    \widetilde{D}& =& (M_4 D + M_3) (M_2 D + M_1)^{-1}.
  \end{array}\right.
\end{equation}
If in addition $(A,B,C,D)$ is balanced, then this realization
$(\widetilde{A},\widetilde{B},\widetilde{C},\widetilde{D})$ is
also balanced.
\end{theorem}
{\em Proof.~}%
A $J$-unitary constant matrix $M$ is in particular $J$-inner and
its inverse is also $J$-inner, so that ${\cal T}_M$ is clearly, by
Prop.\ \ref{LFTonlossless}, a bijection on the set of lossless
functions.
\\
The computation on the realization is just a particular case of
formula (4.85) in \cite{Kimura}. This shows that the McMillan degree does not
increase under the action of ${\cal T}_M$. It cannot decrease
either, because the inverse mapping ${\cal T}_{M^{-1}}={\cal
T}^{-1}_{M}$ can be treated likewise, and also does not increase
the McMillan degree. Hence, the McMillan degree is preserved by
${\cal T}_M$.
\\
Finally, recall that a minimal realization of a lossless function
is balanced if and only if the associated realization matrix is
unitary. Starting from a minimal balanced realization $(A,B,C,D)$
it now can be verified by direct computation that the realization
matrix associated with
$(\widetilde{A},\widetilde{B},\widetilde{C},\widetilde{D})$ is
also unitary and thus the realization is balanced.
\hfill $\Box$%


\section{The tangential Schur algorithm}
\label{tangential}%

In this section we outline the use of the tangential Schur
algorithm for the recursive construction of a parametrization of
the space of stable all-pass systems of fixed finite order. It is
derived from the method of \cite{ABG94,F-O98}, where the
tangential Schur algorithm is used to construct an infinite atlas
of generic overlapping parametrizations for the space of $p \times
p$ inner functions of McMillan degree $n$. The relationship
between these two situations is constituted by the map $R(z)
\mapsto R(z)^{-1}$, or equivalently $R(z) \mapsto R^{\sharp}(z)$
(see also Eqns.\ (\ref{Rsharp}) and (\ref{inverselossless})),
which is used to relate the space of $p \times p$ inner functions
of McMillan degree $n$ to the space of $p \times p$ stable
all-pass systems of McMillan degree $n$.

In the context of inner functions, the tangential Schur algorithm
consists of an iterative procedure by which a given $p \times p$
inner function of McMillan degree $n$ is reduced in $n$ iteration
steps to a $p \times p$ inner function of McMillan degree $0$
(i.e., to a constant unitary matrix). In each iteration step the
McMillan degree of the inner function at hand is reduced by $1$,
by application of a suitable linear fractional transformation
which is chosen to meet a particular interpolation condition, and
which involves an associated $J$-inner matrix function of McMillan
degree $1$ (or elementary $J$-inner factor, see Section \ref{preliminaries}). 
The actual parametrization procedure consists of the
reverse process, by which a chart of inner functions of McMillan
degree $n$ is constructed in $n$ iteration steps, starting from an
initial unitary matrix. The choice of interpolation points and
(normalized) direction vectors may serve to index such a chart,
while the local coordinates correspond to the $n$ Schur vectors.

In the Schur algorithm, the elementary $J$-inner factors involved 
are of the form (\ref{XBlaschke}) and  it suffices to consider 
those  which are analytic on the closed unit disk 
(see the remark at the end of the section).  
Such an elementary $J$-inner factor $\Theta(z)$, having
its pole outside the closed unit disk at $z=1/\overline{w}$, can be
represented as:
\begin{equation}
\label{theta}%
  \Theta(u,v,w,\xi,H)(z) =
  \left( I_{2p} + \left(\frac{b_w(z)}{b_w(\xi)}-1 \right)
  \frac{\left[ \begin{array}{c} u \\ v \end{array} \right]
  \left[ \begin{array}{c} u \\ v \end{array} \right]^*
  J}{(1-\|v\|^2)} \right)\, H,
\end{equation}
where \\
(1) $w \in {\mathbb C}$, $|w|<1$, \\
(2) $u \in {\mathbb C}^{p \times 1}$, $\|u\|=1$, \\
(3) $v \in {\mathbb C}^{p \times 1}$, $\|v\|<1$, \\
(4) $\xi \in {\mathbb C}$, $|\xi|=1$, \\
(5) $H$ is a $2p \times 2p$ constant $J$-unitary matrix.

This representation is obtained from (\ref{XBlaschke}) in which some
normalizing conditions have been imposed:
\begin{itemize} 
\vspace{-0.5cm}
\item[(i)] the vector $x$ has been rewritten  $x = \sigma \left[ \begin{array}{c} u
\\ v \end{array} \right]$, where $\sigma$ is a nonzero scalar, $\|u\|=1$ and
thus $\|v\|<1$, since (\ref{P>0}) must be satisfied.
\item[(ii)] an arbitrary complex number   $\xi$ on the unit circle
has been introduced which allows to write the right $J$-unitary factor 
$H$ as $H=\Theta(\xi)$.
\end{itemize}
\vspace{-0.5cm}
Since  any nonzero scaling of the
vector $x$ does not affect the expression (\ref{XBlaschke}), we have that 
\begin{equation}
\label{thetaXx}
\Theta(u,v,w,\xi,H)(z) =
X_x\left(b_w(z)b_w(\xi)^{-1}\right)H=X_x\left(b_w(z)\right)X_x\left(b_w(\xi)\right)^{-1}H,
\end{equation}
and $\Theta(u,v,w,\xi,H)(z)$ is of the form  (\ref{XBlaschke}) up to a right 
constant $J$-unitary factor.  The constant matrix 
$X_x(b_w(\xi)^{-1}) = X_x(b_w(\xi))^{-1}$ is indeed
$J$-unitary, since $|\xi|=1$ yields $|b_w(\xi)|=1$ (see Lemma \ref{XxYx}). 

 Once $\xi$ is chosen arbitrarily on
the unit circle, the representation (\ref{theta}) 
is unique up to a unimodular complex number, which corresponds to
the remaining freedom in the scaling of the vector $x$. It will
turn out that for several of our purposes, the additional freedom
to choose $\xi$ and $H$ is of rather limited use. In the
literature on the parametrization of inner functions it therefore
often happens that fixed choices are made, in particular $\xi=1$
and $H=I_{2p}$ (see, e.g., \cite{ABG94,F-O98}). However,
especially the role played by $H$ cannot be ignored in
establishing the precise connection between the tangential Schur
algorithm and the state-space approach involving the mappings
${\cal F}_{U,V}$ in Section \ref{connection}. Also, for the
construction of atlases of overlapping parametrizations, the
additional freedom to choose $H$ leads to useful new alternatives,
as will be discussed in Section \ref{overlapping}. Therefore, we
will not fix a choice yet and present our results for general
$\xi$ and $H$.

The following lemma establishes a useful new factorization result
for the class of $J$-unitary matrices $X_x(\alpha)$ with
$x^*Jx>0$.
\begin{lemma}
\label{xcheckfactor}%
Let $x \in {\mathbb{C}}^{2p \times 1}$ be such that $x^*Jx>0$.
Let $x$ be written as $x = \sigma \left[ \begin{array}{c} u \\
v \end{array} \right]$, where $\sigma$ is any nonzero scalar such
that $u, v \in {\mathbb{C}}^{p \times 1}$ satisfy $\|u\|=1$ and
$\|v\|<1$. It then holds that
\begin{equation}
  X_x(\alpha) = H(uv^*) X_{\left[ \begin{array}{c} u \\ 0
  \end{array} \right]}(\alpha) H(uv^*)^{-1},
\end{equation}
where $H(uv^*)$ denotes the Halmos extension of the strictly
contractive matrix $uv^*$.
\end{lemma}
{\em Proof.~}%
Note that $\sigma$, $u$ and $v$ with the required properties
always exist, making $uv^*$ into a strictly contractive matrix.
The Halmos extension $H(uv^*)$ is not difficult to compute
explicitly as the Hermitian $J$-unitary matrix
\begin{equation}
  H(uv^*) = \left[ \begin{array}{cc}
  I_p - (1 - \frac{1}{\sqrt{1-\|v\|^2}})uu^* &
  \frac{1}{\sqrt{1-\|v\|^2}}uv^* \\
  \frac{1}{\sqrt{1-\|v\|^2}}vu^* &
  I_p - (1 - \frac{1}{\sqrt{1-\|v\|^2}})\frac{vv^*}{\|v\|^2}
  \end{array} \right],
\end{equation}
and it satisfies 
\[H(uv^*) \left[\begin{array}{r}u\\-v\end{array}\right]=
\sqrt{1-\|v\|^2} \left[\begin{array}{c}u\\0\end{array}\right].\]
It then follows that $X_x(\alpha) H(uv^*) = H(uv^*) +
{\displaystyle \frac{(\alpha-1)}{(1-\|v\|^2)}} \left[
\begin{array}{c} u \\ v \end{array} \right] \left[
\begin{array}{c} u \\ v \end{array} \right]^* J H(uv^*) = H(uv^*)
+ {\displaystyle \frac{(\alpha-1)}{\sqrt{1-\|v\|^2}}} \left[
\begin{array}{c} u \\ v \end{array} \right] \left[
\begin{array}{c} u \\ 0 \end{array} \right]^* = H(uv^*) +
(\alpha-1) H(uv^*) \left[ \begin{array}{c} u \\ 0 \end{array}
\right] \left[ \begin{array}{c} u \\ 0 \end{array} \right]^* J =
H(uv^*) X_{\left[ \begin{array}{c} u \\ 0 \end{array}
\right]}(\alpha)$, which proves the lemma.
\hfill $\Box$%
\\[5mm]
The following proposition plays a central role in our construction
of balanced parametrizations for discrete-time lossless transfer
functions.
\begin{proposition}
\label{thetafact}%
The $J$-inner matrix function $\Theta(u,v,w,\xi,H)(z)$ can be
factorized as:
\begin{equation}
  \Theta(u,v,w,\xi,H)(z) =
  H(uv^*) S_{u,w}(z) S_{u,w}(b_w(\xi))^{-1} H(uv^*)^{-1} H,
\end{equation}
where $H(uv^*)$ denotes the Halmos extension of the strictly
contractive matrix $uv^*$, and where
\begin{equation}
   S_{u,w}(z) := X_{\left[ \begin{array}{c} u \\ 0 \end{array}
  \right]}(b_w(z)) = \left[ \begin{array}{cc}
  I_p - \left( 1-b_w(z) \right) uu^* & 0 \\ 0 & I_p
  \end{array} \right].
\end{equation}
\end{proposition}
{\em Proof.~}%
This is an immediate consequence of the expression (\ref{thetaXx}) of
$\Theta(u,v,w,\xi,H)(z)$, in combination
with Lemma \ref{xcheckfactor}, in which $\alpha$ is replaced by
the Blaschke factor $b_w(z)$.
\hfill $\Box$%
\\[5mm]
The next proposition establishes that application of the linear
fractional transformation associated with the elementary $J$-inner
function $\Theta(u,v,w,\xi,H)(z)$ to a lossless function of
McMillan degree $n$, results in a lossless function of McMillan
degree $n+1$ which satisfies an interpolation condition involving
only $u$, $v$ and $w$.
\begin{proposition}
Let $\Theta(z)=\Theta(u,v,w,\xi,H)(z)$ be the elementary $J$-inner
function defined in Eqn.\ (\ref{theta}), with $u$, $v$, $w$, $\xi$
and $H$ satisfying the accompanying properties. Let $G(z)$ be a
lossless function of McMillan degree $n$. Then
\begin{equation}
  \widehat{G}(z) = {\cal T}_{\Theta(z)}(G(z)) =
  \widehat{\cal T}_{\Theta^\sharp(z)}(G(z))
\end{equation}
is lossless of McMillan degree $n+1$ and satisfies the interpolation
condition
\begin{equation}
\label{intcond}%
  \widehat{G}(1/\overline w) u = v.
\end{equation}
\end{proposition}
{\em Proof.~}%
Proposition \ref{LFTonlossless} asserts that $\widehat{G}(z)$ is
lossless. From  Proposition \ref{thetafact} we have that
\[
  {\cal T}_{\Theta(z)} = {\cal T}_{H(uv^*)} \circ
  {\cal T}_{S_{u,w}(z)} \circ {\cal T}_{M},
\]
where $M = S_{u,w}(\xi)^{-1} H(uv^*)^{-1} H$. The linear
fractional transformations associated with the $J$-unitary
matrices $M$ and $H(uv^*)$ are both generalized M\"obius
transformations which leave the McMillan degree unchanged, see
Theorem \ref{moebius}. Only the linear fractional transformation
associated with the matrix $S_{u,w}(z)$ can change the McMillan
degree, but it has a simple form: for the lossless function
$R(z):={\cal T}_M(G(z))$ of McMillan degree $n$ it holds that
\[
  \widehat{R}(z) := {\cal T}_{S_{u,w}(z)}(R(z)) =
  R(z) (I_p-(1-b_w^{\sharp}(z)) uu^*).
\]
Because $I_p-(1-b_w^{\sharp}(z)) uu^*$ is the transfer matrix of a
lossless system and, by (\ref{xcheckprop}), $\det{\widehat{R}(z)}
= b_w^{\sharp}(z) \det{R(z)}$, it follows from Lemma
\ref{degreelossless} and from the fact that stability ensures that
no common factors can occur, that $\widehat{R}(z)$ has McMillan
degree $n+1$. Then also $\widehat{G}(z):={\cal
T}_{H(uv^*)}(\widehat{R}(z))$ has McMillan degree $n+1$.

For the matrix $X_x(\alpha)$ it holds that $x^* J X_x(\alpha) =
\alpha x^* J$. Setting $x=\left[ \begin{array}{c} u \\ v
\end{array} \right]$, $\alpha=b_w(z)$ and $z=w$ this results in
the following identity, in view of Eqn.\ (\ref{thetaXx}):
\begin{equation}
\label{thetaintcond0}%
  \left[ \begin{array}{cc} u^* & -v^* \end{array} \right]
  \Theta(w) = 0.
\end{equation}
Application of the operator $^{\sharp}$ yields the equivalent
relationship
\begin{equation}
\label{thetaintcond1}%
  \Theta^{\sharp}(1/\overline{w}) \left[
  \begin{array}{r} u \\ -v \end{array} \right] = 0.
\end{equation}
Note that $\widehat{G}(z) = \widehat{\cal
T}_{\Theta^{\sharp}(z)}(G(z))$ satisfies the relation
\begin{equation}
\label{LFThatmat}%
  \left[ \begin{array}{cc} \widehat{G}(z) & I_p \end{array}
  \right] = (G(z) \Theta^{\sharp}_3(z)+\Theta^{\sharp}_4(z))^{-1}
  \left[ \begin{array}{cc} G(z) & I_p \end{array} \right]
  \Theta^{\sharp}(z).
\end{equation}
The matrix $G(z)\Theta^{\sharp}_3(z)+\Theta^{\sharp}_4(z)$ is
analytic and invertible everywhere outside the open unit disk,
which follows in a manner entirely analogous to the proof of
invertibility of $\Theta_2(z) G(z)+\Theta_1(z)$ at points of
analyticity outside the unit disk, see the proof of Proposition
\ref{LFTonlossless}. Setting $z=1/\overline{w}$ the interpolation
condition (\ref{intcond}) now follows.
\hfill $\Box$%
\\[5mm]
{\em Remark.~}%
Using relation (\ref{thetaintcond0}) and the fact that $\widehat{G}(z) = {\cal
T}_{\Theta(z)}(G(z))$ can be rewritten in the form
\begin{equation}
\label{LFTmat}%
  \left[ \begin{array}{c} I_p \\ \widehat{G} \end{array} \right]
  = \Theta \left[ \begin{array}{c} I_p \\ G \end{array} \right]
  (\Theta_2 G + \Theta_1)^{-1},
\end{equation}
 yield an
alternative form of the interpolation condition (\ref{intcond}):
\begin{equation}
\label{intcondalt}%
  u^* = v^* \widehat{G}(w)
\end{equation}
provided that the matrix $\Theta_2(w)G(w)+\Theta_1(w)$ is
invertible. Equivalently, this requires that $w$ does not show up
as a pole of $\widehat{G}$. Only then the form (\ref{intcondalt})
is valid, while the form (\ref{intcond}) applies always. This
illustrates the usefulness of the identity $\widehat{\cal
T}_{\Theta^{\sharp}} = {\cal T}_{\Theta}$ when dealing with
analyticity. In a similar vein, note that for $w=0$ the function
$\Theta(u,v,w,\xi,H)$ is not proper and does not admit a
state-space realization. 
However, the matrix function $\Theta(u,v,w,\xi,H)^{\sharp}$ does
admit a state-space realization, and formula (4.85) in \cite{Kimura} 
can be used to compute the action of the associated LFT entirely in terms of
realizations.
\\[5mm]
The following proposition addresses the reverse process and
constructs from a lossless matrix function $G(z)$ of McMillan
degree $n+1$, a new lossless matrix function of reduced degree
$n$. To achieve this, it proceeds from an interpolation condition
of the form (\ref{intcond}). Though this result is well-known
(see, e.g. \cite{BGR}, \cite{ABG94}), we find it of interest to give a more
constructive proof of it.
\begin{proposition}
\label{Schurstep}%
Let $\widehat{G}$ be a $p \times p$ lossless function of McMillan
degree $n+1$ which satisfies an interpolation condition of the
form
\begin{equation}
\label{intpcond}%
  \widehat{G}(1/\overline{w}) u = v,
\end{equation}
in which $w \in {\mathbb C}$ is an interpolation point with
$|w|<1$, $u \in {\mathbb C}^{p \times 1}$ is a direction vector
with $\|u\|=1$ and $v \in {\mathbb C}^{p \times 1}$ is a Schur
vector satisfying $\|v\|<1$.
\\
Let $\xi \in {\mathbb C}$ be an arbitrary number of modulus $1$,
and let $H$ be an arbitrary constant $J$-unitary matrix. Then
$\widehat{G}$ admits the representation
\begin{equation}
  \widehat{G} = {\cal T}_{\Theta(u,v,w,\xi,H)}(G),
\end{equation}
for some $p \times p$ lossless function $G$ of McMillan degree
$n$.
\end{proposition}
{\em Proof.~}%
The interpolation condition (\ref{intpcond}) can be rewritten as
follows:
\begin{eqnarray*}
  \widehat{G}(1/\overline{w}) u = v & \Leftrightarrow & 0 =
  \left[ \begin{array}{cc} \widehat{G}(1/\overline{w}) & I_p
  \end{array} \right] \left[ \begin{array}{r} u \\ -v \end{array}
  \right] \\
  & \Leftrightarrow & 0 = \left[ \begin{array}{cc}
  \widehat{G}(1/\overline{w}) & I_p \end{array} \right] H(uv^*)^{-1}
  H(uv^*) \left[ \begin{array}{r} u \\ -v \end{array} \right] \\
  & \Leftrightarrow & 0 = \left[ \begin{array}{cc}
  \widehat{G}(1/\overline{w}) & I_p \end{array} \right] H(uv^*)^{-1}
  \left[ \begin{array}{c} u \\ 0 \end{array} \right] \\
  & \Leftrightarrow & 0 = {\widehat{\cal
  T}}_{H(uv^*)^{-1}}(\widehat{G})(1/\overline{w}) u.
\end{eqnarray*}
The lossless function $\widehat{R} := {\widehat{\cal
T}}_{H(uv^*)^{-1}}(\widehat{G}) = {\cal
T}_{H(uv^*)^{-1}}(\widehat{G})$ has the same McMillan degree $n+1$
as $\widehat{G}$ and satisfies the interpolation condition
\[
  \widehat{R}(1/\overline{w}) u = 0.
\]
This means that $\widehat{R}(z)$ has a zero at $z=1/\overline{w}$,
or equivalently a pole at $z=w$. Now let $B_w(z)$ be defined as
\begin{equation}
  B_w(z) := I_p-(1-b_w^{\sharp}(z))uu^* = (I_p-(1-b_w(z))uu^*)^{-1},
\end{equation}
which is $p \times p$ lossless of McMillan degree $1$, satisfying
$\det B_w(z) = b_w^{\sharp}(z)$ and having its pole at $z=w$. Then
$\widehat{R}(z)$ can be factored as
\begin{equation}
\label{Potfac}%
  \widehat{R}(z) = R(z) B_w(z),
\end{equation}
where by construction $R(z)$ is $p \times p$ lossless of McMillan
degree $n$. (This is a matrix version of the Schwartz lemma on
which the Potapov factorization of lossless and inner functions is
based; see \cite{Pot60}.) The identity (\ref{Potfac}) can be
rewritten as
\[
  \widehat{R} = {\cal T}_{S_{u,w}}(R)
\]
so that finally
\[
  \widehat{G} = {\cal T}_{\Theta(u,v,w,\xi,H)}(G),
\]
where $G={\cal T}_{M^{-1}}(R)$, with $M=S_{u,w}(\xi)^{-1}
H(uv^*)^{-1} H$, is also lossless and of McMillan degree $n$.
\hfill $\Box$%
\\[5mm]
{\em Remark.~}%
In the standard case with $w=0$ the value of $\widehat{G}(\infty)$
corresponds to the direct feedthrough term $\widehat{D}$ of any
state-space realization
$(\widehat{A},\widehat{B},\widehat{C},\widehat{D})$ of
$\widehat{G}(z)$, so that the interpolation condition then takes
the form $\widehat{D}u=v$.
\\[5mm]
The tangential Schur algorithm can now be described. It is used in
Section \ref{overlapping} to construct atlases of generic charts
for the manifold of $p \times p$ lossless systems of degree $n$.
\begin{theorem}[Tangential Schur algorithm]
\label{schuralgo}%
Let $G^{(n)}$ be a $p \times p$ lossless transfer matrix of
McMillan degree $n$. For $k=n,\ldots,1$, let interpolation points
$w_k \in {\mathbb C}$ be given with $|w_k|<1$, let constants
$\xi_k \in {\mathbb C}$ be given with $|\xi_k|=1$, and let
mappings $H_k: (u,v,w,\xi) \mapsto H_k(u,v,w,\xi)$ be given which
assign a constant $2p \times 2p$ $J$-unitary matrix
$H_k(u,v,w,\xi)$ to each quadruple $(u,v,w,\xi) \in {\mathbb C}^p
\times {\mathbb C}^p \times {\mathbb C} \times {\mathbb C}$ with
$\|u\|=1$, $\|v\|<1$, $|w|<1$ and $|\xi|=1$.
\\
Then for $k=n,\ldots,1$, there exist vectors $u_k \in {\mathbb
C}^p$ with $\|u_k\|=1$, such that the vectors $v_k \in {\mathbb
C}^p$ constructed recursively by the following formulas, all have
length $\|v_k\|<1$:
\begin{eqnarray}
\label{vk}%
  & & v_k := G^{(k)}(1/\overline{w}_k) u_k, \\
\label{thetak}%
  & & \Theta_k := \Theta(u_k,v_k,w_k,\xi_k,H_k(u_k,v_k,w_k,\xi_k)), \\
\label{Gk}%
  & & G^{(k-1)} := {\cal T}^{-1}_{\Theta_{k}}(G^{(k)}).
\end{eqnarray}
With such a choice of the unit vectors $u_k$ ($k=n,\ldots,1$) each
of the functions $G^{(k)}$ is lossless of McMillan degree $k$ and
one can write
\begin{equation}
  G^{(n)} = {\cal T}_{\Theta_n}({\cal T}_{\Theta_{n-1}}(\ldots
  ({\cal T}_{\Theta_1}(G^{(0)}))\ldots)) =
  {\cal T}_{\Theta_n \Theta_{n-1}\cdots\Theta_1}(G^{(0)}),
\end{equation}
where $G^{(0)}$ is a constant unitary matrix.
\end{theorem}
{\em Proof.~}%
In view of Proposition \ref{Schurstep}, which is applied
repeatedly, it only remains to show that if $G$ is a lossless
function of degree $n>0$, then for every choice of interpolation
point $w$ with $|w|<1$, a direction vector $u$ with $\|u\|=1$
exists for which the Schur vector $v:=G(1/\overline{w})u$
satisfies $\|v\|<1$. A function theoretic proof of this fact is
given in \cite[Lemma 6]{F-O98}. An alternative state-space
argument runs as follows.
\\
For $|w|<1$, the mapping $z \mapsto b_{-w}(z)$ constitutes an
automorphism of the unit disk. Thus, the function
$\widetilde{G}(z):=G(b_{-w}(z))$ is again lossless (and also of
McMillan degree $n>0$). This shows that
$G(1/\overline{w})=\widetilde{G}(\infty)$ appears as the direct
feedthrough matrix $\widetilde{D}$ of any state-space realization
$(\widetilde{A},\widetilde{B},\widetilde{C},\widetilde{D})$ of the
lossless function $\widetilde{G}$. If such a realization is chosen
to be minimal and balanced, then $\widetilde{D}$ is a proper
submatrix of the associated unitary realization matrix; see
Proposition \ref{unitaryrealization}. Clearly, some column of
$\widetilde{D}$ must have length $<1$, since otherwise all columns
of $\widetilde{B}$ would be zero, contradicting minimality with
order $n>0$. Thus, a suitable vector $u$ can always be found among
the set of standard basis vectors $\{e_1,\ldots,e_p\}$. As a
matter of fact, any normalized basis will do.
\hfill $\Box$%
\\[5mm]
{\em Remark.~}%
Similar results can be proved in which the interpolation condition
(\ref{intcond}), which involves a relationship among the {\em
columns} of $\widehat{G}$, is replaced by an interpolation
condition involving its {\em rows}. This can be done by applying
Proposition \ref{Schurstep} to the lossless function
$\widehat{G}^*(z)$ and using part (c) of Lemma
\ref{LFTproperties}. The interpolation condition then takes the
form
\begin{equation}
\label{leftintcond}%
  u^* \widehat{G}(1/w) = v^*,
\end{equation}
and should not be confused with the alternative form
(\ref{intcondalt}) for the interpolation condition
(\ref{intcond}). As a counterpart to the result of Proposition
\ref{Schurstep}, one also obtains an alternative linear fractional
representation
\[
  \widehat{G} = {\cal T}_{\Theta(u,v,w,\xi,H)^o}(\widetilde{G})
\]
for some $p \times p$ lossless function $\widetilde{G}$ of degree
$n$. Observe that the $J$-inner function $\Theta(u,v,w,\xi,H)^o$
is now analytic outside the open unit disk. Thus, the theory of
this section admits a dual theory in which elementary $J$-inner
factors analytic outside the disk and interpolation conditions on
the rows are used.
\\
Note that the tangential Schur algorithm described above proceeds
from a set of $n$ interpolation conditions on the columns of a
sequence of intermediate lossless functions. Obviously, a set of
$n$ interpolation conditions on the rows of a sequence of
intermediate lossless functions can also be imposed. But since
each iteration step involves a single interpolation condition
only, {\em mixed} situations can also be treated in which $n$
interpolation conditions are prescribed involving columns and
rows.


\section{Connection between the classes of mappings
         ${\cal T}_{\Theta(z)}$ and ${\cal F}_{U,V}$}
\label{connection}%

In the scalar case, mappings of the form ${\cal F}_{U,V}$ with
$U=I_2$ have been used in \cite{HP2000} to recursively construct a
balanced canonical form for the space of discrete-time lossless
systems of finite McMillan degree. The parameters that occur in
this recursion have the interpretation of Schur parameters,
corresponding to the situation with interpolation points $w$ at
zero. With this connection in mind, it is the purpose of the
present section to clarify the relationship between the two
classes of mappings ${\cal T}_{\Theta(z)}$ and ${\cal F}_{U,V}$,
with $\Theta(z)$ a $J$-inner matrix function of McMillan degree
$1$ and $U$ and $V$ unitary. We are interested in investigating
the possibilities for representing a mapping ${\cal
T}_{\Theta(z)}$ in terms of a corresponding mapping ${\cal
F}_{U,V}$. This will give us balanced state-space parametrizations
directly in terms of the set of parameters $u$, $v$, $w$, $\xi$
and $H$ used in the tangential Schur algorithm. Moreover, unitary
matrices will be involved in the computation of these
realizations, and these are known to be numerically
well-conditioned.

Let $\cL_n^p$ denote the set of rational $p
\times p$ lossless transfer matrices of McMillan degree $n$ and
let $\cL^p=\bigcup_{n=0}^{\infty} \cL_n^p$.
\begin{theorem}
\label{FUVLFT}%
Let $U$ and $V$ be unitary $(p+1) \times (p+1)$ matrices,
block-partitioned as
\begin{equation}
\label{UVpartitions}%
  U = \left[ \begin{array}{cc} \alpha_U & M_U \\ k_U & \beta_U^*
  \end{array} \right],
  \hspace{5mm}
  V = \left[ \begin{array}{cc} \alpha_V & M_V \\ k_V & \beta_V^*
  \end{array} \right],
\end{equation}
with $k_U$ and $k_V$ scalar and the other blocks of compatible sizes.
\\
 The mapping ${\cal F}_{U,V}: \cL^p \rightarrow \cL^p$ can be
expressed as a linear fractional transformation ${\cal T}_{\Phi}:
\cL^p \rightarrow \cL^p$ for some $J$-inner function $\Phi$ if and
only if $(k_U,k_V) \neq (0,0)$.
\\
(i) If $|k_U|=|k_V|=1$, then $M_U$ and $M_V$ are unitary matrices
and ${\cal F}_{U,V}$ is the generalized M\"{o}bius transform
${\cal F}_{U,V}(G) = M_V G M_U^*$ which can be represented as
${\cal T}_{\Phi}(G)$ with a constant $J$-unitary matrix $\Phi$.
\\
(ii) If $(k_U,k_V) \neq (0,0)$ and $|k_U|<1$ or $|k_V|<1$, then
${\cal F}_{U,V}$ can be represented as ${\cal T}_{\Phi}$ with an
elementary $J$-inner function $\Phi$ having a pole at $k_V/k_U$,
located at infinity for $k_U=0$.
\end{theorem}
{\em Proof.~}%
Because $U$ and $V$ are unitary, according to Proposition
\ref{FUVonlossless} the mapping ${\cal F}_{U,V}$ maps $\cL^p$ into
$\cL^p$. If $\Phi$ is $J$-inner, then the mapping ${\cal
T}_{\Phi}$ also maps $\cL^p$ into $\cL^p$ but in addition it is
bijective, see Lemma \ref{That} and Proposition
\ref{LFTonlossless}. Thus, a necessary requirement for ${\cal
F}_{U,V}$ to be representable as a linear fractional
transformation ${\cal T}_{\Phi}$ with a $J$-inner function $\Phi$,
is injectivity.
\\
First consider the scalar case $p=1$ with $(k_U,k_V)=(0,0)$. Then
$U$ and $V$ are $2 \times 2$ diagonal matrices and all lossless
$G$ are mapped to the unimodular constant $\alpha_V\alpha_U^*$.
Therefore the mapping ${\cal F}_{U,V}$ obviously is not injective.
\\
In the multivariable case $p>1$ with $(k_U,k_V)=(0,0)$, it also
holds that ${\cal F}_{U,V}$ is not injective and thus cannot be
represented as an LFT associated with a $J$-inner function.
Indeed, let $P(\xi)=I_p+(\xi-1)\beta_U\beta_U^*$, and
$Q(\xi)=I_p+(\xi-1)\beta_V\beta_V^*$ where $\xi \neq 1$ is a
unimodular number. Then $P(\xi)$ and $Q(\xi)$ are unitary and we
have
\[
  U \left[ \begin{array}{cc} 1 & 0 \\ 0 & P(\xi) \end{array}
  \right] = \left[ \begin{array}{cc} I_p & 0 \\ 0 & \xi
  \end{array} \right] U, \hspace{5mm}
  V \left[ \begin{array}{cc} 1 & 0 \\ 0 & Q(\xi) \end{array}
  \right] = \left[ \begin{array}{cc} I_p & 0 \\ 0 & \xi
  \end{array} \right] V,
\]
so that by Lemma \ref{3properties} (i), (iii), for every lossless
function $G$,
\[
  {\cal F}_{U,V}(Q(\xi)G(z)P(\xi)^*) =
  {\cal F}_{\left[ \begin{array}{cc} I_p & 0 \\ 0 & \xi
  \end{array} \right] U, \left[ \begin{array}{cc} I_p & 0 \\
  0 & \xi \end{array} \right] V} \left( G(z) \right) =
  {\cal F}_{U,V}(G(z)).
\]
Therefore injectivity of ${\cal F}_{U,V}$ requires that $Q(\xi) G
P(\xi)^* = G$ for every lossless function $G$. However if this is
the case, then $\beta_V \beta_V^* G = G \beta_U \beta_U^*$ for all
lossless $G$. Choosing $G=I_p$ this implies $\beta_V \beta_V^* =
\beta_U \beta_U^*$. Choosing $G=I_p - 2 \gamma \gamma^*$, where
$\gamma \in {\mathbb C}^p$ is a vector of length one which is
independent of $\beta_V$ and not orthogonal to $\beta_V$, it
follows that $\beta_V \beta_V^* \gamma \gamma^* = \gamma \gamma^*
\beta_U \beta_U^*$, which yields a contradiction. Therefore ${\cal
F}_{U,V}$ is not injective.
\\
Next, assume that $(k_U,k_V) \neq (0,0)$. Let
\[
  J_1 = \left[ \begin{array}{cc} I_p & 0 \\ 0 & 0 \end{array}
  \right] \hspace{5mm} \mbox{\rm and} \hspace{5mm}
  J_2 = \left[ \begin{array}{cc} 0 & 0 \\ 0 & 1 \end{array}
  \right]
\]
and consider the linear fractional transformation ${\cal L}$
defined by
\[
  {\cal L}(H(z)) = (J_1 H(z) + J_2) (J_2 H(z) + J_1)^{-1}.
\]
This linear fractional transformation is well-defined for every
$H(z)$ such that, in the block decomposition of $H(z)$
\[
  H(z) = \left[ \begin{array}{cc} H_1(z) & H_2(z) \\
  H_3(z) & H_4(z) \end{array} \right],
\]
with $H_4(z)$ scalar, $H_4(z)$ is not identically zero. It can be
computed as
\[
  {\cal L}(H) =
  \left[ \begin{array}{cc} H_1- H_2 H_4^{-1}H_3 & H_2 H_4^{-1} \\
  - H_4^{-1} H_3 & H_4^{-1} \end{array} \right].
\]
The mapping ${\cal F}_{U,V}$ can now be written using this
transformation as:
\[
  {\cal F}_{U,V}(G(z)) = \left[ \begin{array}{ll} I_p & 0
  \end{array} \right] {\cal L} \left( V \left[ \begin{array}{cc}
  1 & 0 \\ 0 & G(z) \end{array} \right] U^* -z \, \left[
  \begin{array}{cc} 0 & 0 \\ 0 & 1 \end{array} \right] \right)
  \left[ \begin{array}{l} I_p \\ 0 \end{array} \right].
\]
Now,
\begin{eqnarray*}
 & & {\cal L} \left( V \left[ \begin{array}{cc} 1 & 0 \\ 0 & G(z)
 \end{array} \right] U^* - z \, \left[ \begin{array}{cc} 0 & 0 \\
 0 & 1 \end{array} \right] \right)= \\
 & & = \left[ J_1 \left( V \left[ \begin{array}{cc} 1 & 0 \\ 0 &
 G(z) \end{array} \right] - z \, \left[ \begin{array}{cc} 0 & 0 \\
 0 & 1 \end{array} \right] U \right) + J_2 U \right]
 \left[ J_2 \left( V \left[ \begin{array}{cc} 1 & 0 \\ 0 & G(z)
 \end{array} \right] - z \, \left[ \begin{array}{cc} 0 & 0 \\
 0 & 1 \end{array} \right] U \right) + J_1 U \right]^{-1}= \\
 & & = \left[ \begin{array}{cc} \alpha_V & M_V G \\ k_U &
 \beta_U^* \end{array} \right] \left[ \begin{array}{cc}
 \alpha_U & M_U \\ k_V - z \, k_U & \beta_V^* G - z \, \beta_U^*
 \end{array} \right]^{-1}.
\end{eqnarray*}
Using the following expression for the inverse of a partitioned
matrix
\[
  \left[\begin{array}{cc} M_1 & M_2 \\ M_3 & M_4 \end{array}
  \right] = \left[ \begin{array}{cc} 0 & M_3^{-1} \\ 0 & 0
  \end{array} \right] + \left[ \begin{array}{c} -M_3^{-1}M_4 \\
  I_p \end{array} \right] (M_2-M_1M_3^{-1}M_4)^{-1} \left[
  \begin{array}{cc} I_p & -M_1M_3^{-1} \end{array} \right],
\]
in which $M_2$ and $M_3$ are square submatrices, and $M_3^{-1}$ is
assumed to exist (which holds by assumption in our case), we
obtain
\begin{eqnarray*}
  & & {\cal F}_{U,V}(G) = \left[ \begin{array}{cc} \alpha_V & M_V
  \, G \end{array} \right] \left[ \begin{array}{c}
  -\frac{\beta_V^* \, G - z \, \beta_U^*}{k_V - z \, k_U} \\ I_p
  \end{array} \right] \left[ M_U -\frac{\alpha_U \, (\beta_V^*
  \, G - z \, \beta_U^*)}{k_V - z \, k_U} \right]^{-1} = \\
  & & = \left[ M_V \, G - \frac{\alpha_V \, (\beta_V^* \, G - z
  \, \beta_U^*)}{k_V - z \, k_U} \right] \left[ M_U -
  \frac{\alpha_U \, (\beta_V^* \, G - z\, \beta_U^*)}{k_V - z \,
  k_U} \right]^{-1} = \\
  & & = \left[ \left( M_V \, -\frac{\alpha_V \, \beta_V^*}{k_V - z
  \, k_U} \right) \, G + \frac{z \, \alpha_V \, \beta_U^*}{k_V - z
  \, k_U} \right] \left[-\frac{\alpha_U \, \beta_V^*}{k_V -z \,
  k_U} \, G + \left( M_U + \frac{z \, \alpha_U \, \beta_U^*}{k_V -
  z \, k_U} \right) \right]^{-1}.
\end{eqnarray*}
And finally, ${\cal F}_{U,V}(G) = {\cal T}_{\Phi}(G)$ for $\Phi$
given by
\begin{equation}
\label{PhiUV}%
  \Phi(z) = \left[ \begin{array}{cc}
  M_U + \displaystyle \frac{z \alpha_U \beta_U^*}{k_V - z k_U} &
  \displaystyle -\frac{\alpha_U \beta_V^*}{k_V - z k_U} \\
  \displaystyle \frac{z \alpha_V \beta_U^*}{k_V - z k_U} &
  M_V \displaystyle -\frac{\alpha_V \beta_V^*}{k_V - z k_U}
  \end{array} \right].
\end{equation}
Exploiting unitarity of $U$ and $V$, it then is straightforward to
establish that
\[
  J - \Phi(z) J \Phi(\lambda)^* =
  \frac{(1-\overline{\lambda} \, z)}{(k_V - z \, k_U)
  (\overline{k}_V - \overline{\lambda} \, \overline{k}_U)}
  \left[ \begin{array}{c} \alpha_U \\ \alpha_V \end{array} \right]
  \left[ \begin{array}{c} \alpha_U \\ \alpha_V \end{array}
  \right]^*,
\]
from which it follows that $\Phi(z)$ is $J$-inner. The expression
(\ref{PhiUV}) for $\Phi(z)$ is of McMillan degree $1$ having a
pole at $k_V/k_U$ (located at infinity if $k_U=0$), except if
$|k_U|=|k_V|=1$, in which case $\Phi(z)$ is of McMillan degree $0$
since $\alpha_U=\alpha_V=0$ and $\Phi(z) = \left[\begin{array}{cc}
M_U & 0 \\ 0 & M_V \end{array} \right]$ becomes constant.
$\mbox{}$\hfill $\Box$%
\\[5mm]
The previous theorem makes clear under which conditions a mapping
${\cal F}_{U,V}$ with unitary matrices $U$ and $V$ admits a
representation as an LFT associated with an elementary $J$-inner
factor $\Phi$.
Depending on the location of the pole $k_V/k_U$, this 
factor takes the different forms described in Theorem \ref{Jinnerfactors}.
In the cases where such a representation exists the
proof is constructive, yielding the possible expression
(\ref{PhiUV}) for $\Phi$. Using Lemma \ref{LFTlosslessequiv}
\emph{all} the possibilities for $\Phi$ can be computed.
\begin{lemma}
\label{uniqueJelementary}%
Let $\Phi(z)$ and $\Psi(z)$ be elementary $J$-inner factors for
which ${\cal T}_{\Phi}$ and ${\cal T}_{\Psi}$ agree on $\cL^p$. If
$p>1$ then there exists a unimodular constant $\rho$ such that
\[
  \Phi(z) = \rho \Psi(z).
\]
If $p=1$, then either there exists a unimodular constant $\rho$
such that $\Phi(z) = \rho \Psi(z)$, or $\Phi(z)$ has a pole at
$z=w$ with $|w| \neq 1$, $\Psi(z)$ has a pole at
$z=1/\overline{w}$ and there exists a unimodular constant $\rho$
such that
\[
  \Phi(z) = \rho b_w^\sharp(z) \Psi(z).
\]
\end{lemma}
{\em Proof.~}%
By Lemma \ref{LFTlosslessequiv},
\begin{equation}
\label{equiv-Jinner}%
  \Phi(z) = b(z) \Psi(z),
\end{equation}
for some scalar function $b(z)$, and since $\Phi$ and $\Psi$ are
$J$-inner, $b(z)$ must be a Blaschke factor. Computing the
determinants gives
\[
  \det \Phi(z) = b(z)^{2p} \det \Psi(z).
\]
But the determinant of an elementary $J$-inner function is either
equal to a unimodular number (if the pole is located on the unit
circle) or to a Blaschke factor, having McMillan degree $1$. If
$p>1$ it therefore follows that that $b(z)$ must be a unimodular
constant. If $p=1$ an alternative possibility occurs if $\Phi$ has
a pole at $z=w$ with $|w| \neq 1$. In that case $\Psi$ may have a
pole at $z=1/\overline{w}$ and $b(z)$ may also be chosen equal to
$\rho b_w^{\sharp}(z)$ for some complex number $\rho$ of modulus
one. If $\Psi$ is of the form
\[
  \Psi(z) = \left( I_2 - (1-b_w(z)) x (x^*Jx)^{-1} x^*J \right) H,
\]
then
\[
  \Phi(z) = \rho \left( I_2 - (1-b_w{^\sharp}(z)) y (y^*Jy)^{-1}
  y^*J \right) H,
\]
where $y$ is a solution to $y^*Jx=0$.
\hfill $\Box$%
\\[5mm]
As a consequence of Theorem \ref{FUVLFT} and Lemma
\ref{uniqueJelementary} we have the following result.
\begin{corollary}
\label{fuvaslft}%
Let $U$ and $V$ be unitary matrices of size $(p+1) \times (p+1)$,
partitioned as in (\ref{UVpartitions}). \\
If $p>1$, then the mapping ${\cal F}_{U,V}$ can be represented in
the form of a mapping ${\cal T}_{\Theta(u,v,w,\xi,H)}$ if and only
if $|k_U|<|k_V|$. Likewise it can be represented in the form of a
mapping ${\cal T}_{\Theta(u,v,w,\xi,H)^o}$ if and only if
$|k_U|>|k_V|$.
\\
If $p=1$, a representation of ${\cal F}_{U,V}$ in the form of a
mapping ${\cal T}_{\Theta(u,v,w,\xi,H)}$ exists if and only if a
representation in the form of a mapping ${\cal
T}_{\Theta(u,v,w,\xi,H)^o}$ exists, which holds if and only if
$|k_U| \neq |k_V|$.
\end{corollary}
{\em Proof.~}%
The statements for the case $p>1$ follow from the earlier results
as indicated. For $p=1$, note that ${\cal F}_{U,V} = {\cal
F}_{\overline{V},\overline{U}}$ which follows from Lemma
\ref{3properties} (ii) upon complex conjugation and noting that
transposition has no effect on a matrix of size $1 \times 1$.
\hfill $\Box$%
\\[5mm]
To complete the picture it remains to characterize the elementary
$J$-inner factor $\Phi$ analytic inside the disk in case
$|k_U|<|k_V|$. As we have seen, if a mapping ${\cal F}_{U,V}$ is
given, the expression (\ref{PhiUV}) together with Lemma
\ref{uniqueJelementary} provide all the possible alternatives for
$\Phi$. Conversely, if an elementary $J$-inner factor $\Phi$
analytic inside the disk is given, one may wonder whether ${\cal
T}_{\Phi}$ admits a representation as a mapping of the form ${\cal
F}_{U,V}$. Here it will prove to be essential to introduce the
$J$-inner matrix function $\widehat{\Theta}(u,v,w)(z)$ defined by
\begin{equation}
  \widehat{\Theta}(u,v,w)(z) =
  H(uv^*) S_{u,w}(z) H(\overline{w}uv^*).
\end{equation}
The main result of this section can now be stated as follows.
\begin{theorem}
\label{generalcase}%
Let $\Theta(z)$ be an elementary $J$-inner factor analytic inside
the unit disk. If ${\cal T}_{\Theta(z)}$ coincides with a mapping
${\cal F}_{U,V}$ with $U$ and $V$ unitary, then $\Theta(z)$ can be
written as
\begin{equation}
  \Theta(z) = \widehat{\Theta}(u,v,w)(z)
  \left[ \begin{array}{cc} P & 0 \\ 0 & Q \end{array} \right],
\end{equation}
for some $u, v \in {\mathbb C}^p$ with $\|u\|=1$, $\|v\|<1$, some
$w \in {\mathbb C}$ with $|w|<1$, and some $p \times p$ unitary
matrices $P$ and $Q$.
\\
In that case one can take $U = \widehat{U} \left[
\begin{array}{cc} 1 & 0 \\ 0 & P \end{array} \right]$ and $V =
\widehat{V} \left[ \begin{array}{cc} 1 & 0 \\ 0 & Q \end{array}
\right]$, where
\begin{eqnarray}
\label{uhat}%
  & & \widehat{U} = \left[ \begin{array}{cc}
      \frac{\sqrt{1-|w|^2}}{\sqrt{1-|w|^2\|v\|^2}} u &
      I_p - (1+\frac{w\sqrt{1-\|v\|^2}}{\sqrt{1-|w|^2\|v\|^2}})uu^* \\
      \frac{\overline{w}\sqrt{1-\|v\|^2}}{\sqrt{1-|w|^2\|v\|^2}} &
      \frac{\sqrt{1-|w|^2}}{\sqrt{1-|w|^2\|v\|^2}} u^*
      \end{array} \right], \\
  & & \nonumber \\
\label{vhat}%
  & & \widehat{V} = \left[ \begin{array}{cc}
      \frac{\sqrt{1-|w|^2}}{\sqrt{1-|w|^2\|v\|^2}} v &
      I_p - (1-\frac{\sqrt{1-\|v\|^2}}{\sqrt{1-|w|^2\|v\|^2}})
            \frac{vv^*}{\|v\|^2} \\
      \frac{\sqrt{1-\|v\|^2}}{\sqrt{1-|w|^2\|v\|^2}} &
      -\frac{\sqrt{1-|w|^2}}{\sqrt{1-|w|^2\|v\|^2}} v^*
      \end{array} \right].
\end{eqnarray}
\end{theorem}
{\em Proof.~}%
Theorem \ref{FUVLFT} ensures that a mapping ${\cal F}_{U,V}$ can
be represented by a linear fractional transformation for the
$J$-inner elementary factor $\Phi(z)$ given by (\ref{PhiUV}). If
$p>1$, it follows from Lemma \ref{uniqueJelementary} that
$|k_U|<|k_V|$. If $p=1$, since  the mappings ${\cal F}_{U,V}$ and
${\cal F}_{\overline{V},\overline{U}}$ coincide, we still may
assume $|k_U|<|k_V|$, and in both cases we have that
\[
  \Phi(z) = \rho \, \Theta(z),
\]
for some unimodular number $\rho$.
\\
Since $|k_U|<|k_V|$, $k_V \neq 0$ and $k_U \neq 1$, so that
$\|\alpha_U\| \neq 0$. Let
\begin{equation}
\label{wuv}%
  w = \frac{\overline{k}_U}{\overline{k}_V}, \hspace{5mm}
  u = \frac{\alpha_U}{\|\alpha_U\|}, \hspace{5mm}
  v = \frac{\alpha_V}{\|\alpha_U\|},
\end{equation}
which satisfy $|w|<1$, $\|u\|=1$ and $\|v\|<1$. The function
$\Theta(z)$ is, by Proposition \ref{thetafact}, of the form
\[
  \Theta(z) = H(uv^*) S_{u,w}(z) H(E) \left[ \begin{array}{cc}
  P & 0 \\ 0 & Q \end{array} \right],
\]
for some strictly contractive matrix $E$ and some unitary matrices
$P$ and $Q$, depending on the matrices $U$ and $V$. To see this,
note that if $\Theta(z)$ is an elementary $J$-inner factor written
in the form $\Theta(u,v,w,\xi,H)(z)$ having its only zero at
$z=w$, then the matrix $\Theta(w)$ has a nontrivial left kernel of
dimension one, which is characterized by Eqn.\
(\ref{thetaintcond0}). For the matrix function $\Phi(z)$ given by
(\ref{PhiUV}) it is easily shown that $\left[ \begin{array}{cc}
u^* & -v^* \end{array} \right] \Phi(w) = 0$ for $w$, $u$ and $v$
defined by Eqn.\ (\ref{wuv}).
\\
Note that, the matrices $U$ and $V$ being unitary,
\begin{eqnarray*}
   & & \|\alpha_U\|^2 + |w|^2|k_V|^2 = 1 \\
   & & \|\alpha_U\|^2\|v\|^2 + |k_V|^2 = 1,
\end{eqnarray*}
so that
\begin{equation}
\label{UVcolumn1}%
  \|\alpha_U\| = \frac{\sqrt{1-|w|^2}}{\sqrt{1-|w|^2\|v\|^2}},
  \hspace{5mm}
  |k_V| = \frac{\sqrt{1-\|v\|^2}}{\sqrt{1-|w|^2\|v\|^2}}.
\end{equation}
Now, the quickest way to proceed is to compare the realizations of
$\Phi^{\sharp}$ and $\Theta^{\sharp}$ (note that $\Theta$ will
fail to have a realization if $w=0$). Using (\ref{PhiUV}) we get
\[
  \Phi^{\sharp}(z)=
  \left[ \begin{array}{cc} M_U^* & 0 \\
  -\frac{\beta_V \alpha_U^*}{\overline{k}_V} & M_V^* -
  \frac{\beta_V \alpha_V^*}{\overline{k}_V} \end{array} \right] +
  \left[ \begin{array}{c} \beta_U \\
  -\frac{\overline{k}_U}{\overline{k}_V} \beta_V \end{array}
  \right] (z-\frac{\overline{k}_U}{\overline{k}_V})^{-1}
  \left[ \begin{array}{cc} \frac{\alpha_U^*}{\overline{k}_V} &
  \frac{\alpha_V^*}{\overline{k}_V} \end{array} \right],
\]
while a state-space realization of $S_{u,w}^{\sharp}(z)$ is easily
computed to yield that
\[
  \Theta^{\sharp}(z) = {\cal D} + {\cal C}(z-w)^{-1}{\cal B},
\]
where
\[
  {\cal D} = \left[ \begin{array}{cc} P^* & 0 \\ 0 & Q^*
  \end{array} \right] H(E)
  \left[ \begin{array}{cc} I_p-(1+\overline{w}) uu^* & 0 \\
  0 & I_p \end{array} \right] H(uv^*),
\]
\[
  {\cal C} = \left[ \begin{array}{cc} P^* & 0 \\ 0 & Q^*
  \end{array} \right] H(E)
  \left[ \begin{array}{c} \sqrt{1-|w|^2}u \\ 0 \end{array}
  \right] \hspace{5mm} \mbox{\rm and} \hspace{5mm}
  {\cal B} = \left[\begin{array}{cc} \sqrt{1-|w|^2}u^* & 0
  \end{array} \right] H(uv^*).
\]
Comparing the (1,2)-block-entries of the expressions for the
direct feedthrough yields, using (\ref{Halmos}),
\[
   P^*(I_p-EE^*)^{-1/2} (E-\overline{w} u v^*) (I_p-vv^*)^{-1/2} = 0
\]
and thus $E=\overline{w} u v^*$ or equivalently
\[
  \Theta(z) = \widehat{\Theta}(u,v,w)(z)
  \left[ \begin{array}{cc} P & 0 \\ 0 & Q \end{array} \right],
\]
as claimed. Replacing $E$ by its value gives
\[
  {\cal D} = \left[ \begin{array}{cc} P^* & 0 \\ 0 & Q^*
  \end{array} \right] \left[ \begin{array}{cc}
  I_p - \left( 1 + \frac{\overline{w} \sqrt{1-\|v\|^2}}
  {\sqrt{1-|w|^2\|v\|^2}} \right) u u^* & 0 \\
  \frac{1-|w|^2}{\sqrt{1-\|v\|^2}\sqrt{1-|w|^2\|v\|^2}} v u^* &
  I_p - \left( 1 - \frac{\sqrt{1-|w|^2\|v\|^2}}
  {\sqrt{1-\|v\|^2}} \right) \frac{vv^*}{\|v\|^2} \end{array}
  \right],
\]
\[
  {\cal C} = \frac{\sqrt{1-|w|^2}}{\sqrt{1-|w|^2\|v\|^2}}
  \left[ \begin{array}{cc} P^* & 0 \\ 0 & Q^* \end{array} \right]
  \left[ \begin{array}{c} u \\ w v \end{array} \right]
  \hspace{5mm} \mbox{\rm and} \hspace{5mm}
  {\cal B} = \frac{\sqrt{1-|w|^2}}{\sqrt{1-\|v\|^2}}
  \left[ \begin{array}{cc} u^* & v^* \end{array} \right].
\]
The matrix $M_U$ is therefore given by
\[
  M_U = \rho \left( I_p - \left( 1 + \frac{w \sqrt{1-\|v\|^2}}
  {\sqrt{1-|w|^2\|v\|^2}} \right) uu^* \right) P.
\]

From (\ref{wuv}) and (\ref{UVcolumn1}) we have that
\[
  \left[ \begin{array}{cc} \frac{\alpha_U^*}{\overline{k}_V} &
  \frac{\alpha_V^*}{\overline{k}_V} \end{array} \right] =
  \frac{|k_V|}{\overline{k}_V} \ {\cal B},
\]
which determines the transformation between the two realizations, so that
\[
  \left[ \begin{array}{cc} \beta_U^* & -\frac{k_U}{k_V} \beta_V^*
  \end{array} \right] =
  \rho \ \frac{\overline{k}_V}{|k_V|} \ {\cal C}^* = \rho \
  \frac{\overline{k}_V}{|k_V|} \
  \frac{\sqrt{1-|w|^2}}{\sqrt{1-|w|^2\|v\|^2}}
  \left[ \begin{array}{cc} u^* & \overline{w} \, v^* \end{array}
  \right] \left[ \begin{array}{cc} P & 0 \\ 0 & Q \end{array}
  \right].
\]
Finally it is obtained that the matrices $U$ and $V$ are of the
form
\[
  U = \left[\begin{array}{cc} I_p & 0 \\ 0 &
  \frac{\overline{k}_V}{|k_V|} \end{array} \right] \widehat{U}
  \left[ \begin{array}{cc} 1 & 0 \\ 0 & \rho P \end{array}
  \right] \hspace{5mm} \mbox{\rm and} \hspace{5mm}
  V = \left[\begin{array}{cc} I_p & 0 \\ 0 &
  \frac{\overline{k}_V}{|k_V|} \end{array} \right] \widehat{V}
  \left[ \begin{array}{cc} 1 & 0 \\ 0 & \rho Q \end{array}
  \right],
\]
with $\hat{U}$ and $\hat{V}$ specified by (\ref{uhat}) and
(\ref{vhat}). Since a pre-multiplication by $\left[
\begin{array}{cc} I_p & 0 \\ 0 & \frac{\overline{k}_V}{|k_V|}
\end{array} \right]$ and a post-multiplication by
$\left[\begin{array}{cc} 1 & 0 \\ 0 & \rho I_p \end{array}
\right]$ do not change the mapping ${\cal F}_{U,V}$ according to
Lemma \ref{3properties}, the values of $U$ and $V$ given in the
theorem are obtained.
\hfill $\Box$%


\section{Overlapping canonical forms for lossless and stable systems}
\label{overlapping}%

In the tangential Schur algorithm as presented in Theorem
\ref{schuralgo}, a lossless function $G=G^{(n)}$ of McMillan
degree $n$ is broken down to a constant unitary matrix $G^{(0)}$
using a sequence of interpolation points $w_k$, direction vectors
$u_k$, Schur vectors $v_k$, unimodular constants $\xi_k$ and
constant $J-$unitary matrices $H_k(u_k,v_k,w_k,\xi_k)$, for
$k=1,\ldots,n$. The set of values for $w_k$, $u_k$, $\xi_k$ and
the set of mappings $H_k$ at the $n$ recursion steps can serve to
index a generic chart for the differentiable manifold of lossless
functions of degree $n$, if the mappings $H_k$ are sufficiently
smooth. The Schur parameter vectors $v_k$ together with the
unitary matrix $G^{(0)}$ that we finally reach, then provide the
local coordinates for this chart. An infinite atlas of overlapping
generic charts is obtained by varying the choices for the $u_k$,
$w_k$, $\xi_k$ and $H_k$. Sub-atlases can be extracted from it, as
the one built in \cite{ABG94} which employs elementary $J$-inner
matrices of the form $\Theta(u_k,v_k,w_k,1,I_p)$. Such a choice
implies that $G^{(0)}$ is the value of the lossless function $G$
at $1$. We may even take a finite atlas by fixing all the
interpolation points at zero and letting the values of the $u_k$
vary among a canonical basis of $\CC^p$ (cf.\ the proof of Theorem
\ref{schuralgo}). However it is interesting to keep in mind the
richness of the possibilities, including the one to mix
interpolation conditions on the columns (\ref{intcond}) with
interpolation conditions on the rows (\ref{leftintcond}). This
could be of interest in order to describe lossless transfer
functions with a particular structure.

Now let constants $\xi_k$ be chosen with $|\xi_k|=1$ for all
$k=1,\ldots,n$, and let the mappings 
\[H_k(u,v,w,\xi)= H(uv^*)S_{u,w}(\xi)H(\bar w uv^*)\]
 be chosen so that
$\Theta(u,v,w,\xi_k,H_k(u,v,w,\xi_k)) = \widehat{\Theta}(u,v,w)$
for all $u \in {\mathbb C}^p$, $\|u\|=1$, $v \in {\mathbb C}^p$,
$\|v\|<1$, $w \in {\mathbb C}$, $|w|<1$. Then one obtains an atlas
for the manifold of lossless systems that can be described
directly in terms of state-space realizations using the ${\cal
F}_{U,V}$ mappings.

Let us now describe such an atlas for the manifold $\cL_n^p$ of $p
\times p$ lossless functions of McMillan degree $n$ in more
detail. Choose $\xi_k=1$ for all $k=1,\ldots,n$. For a set of
interpolation points $w_1,w_2,\ldots,w_n$ in the open unit disk, a
set of direction vectors $u_1,u_2,\ldots,u_n \in \CC^p$ of length
one, and a coordinate chart $(\cW,\vartheta)$ for the set $\cU_p$
of unitary $p \times p$ matrices, we define a chart
$(\cV,\varphi)$ by its domain
\[
  \cV{(w_1,w_2,\ldots,w_n,u_1,u_2,\ldots,u_n,\cW)} =
  \{ G \in \cL_n^p \, | \, \|G^{(k)}(1/\overline{w}_k)u_k\|<1,
  G^{(0)} \in \cW \},
\]
and its coordinate map
\[
  \begin{array}{cccl}
  \varphi : & G & \to & (v_1, v_2, \ldots, v_n, \vartheta(G^{(0)}))
  \end{array}
\]
where, for $k=1,\ldots,n$, the lossless functions $G^{(k)}$ and
the Schur parameter vectors $v_k$ are recursively defined in
Theorem \ref{schuralgo}, in which we choose
\[
  \Theta_k = \widehat{\Theta}(u_k,v_k,w_k).
\]
\begin{theorem}
The family $(\cV,\varphi)$ defines a $C^{\infty}$-atlas on
$\cL_n^p$.
\end{theorem}
The proof is left to the reader; see also \cite{ABG94} and
\cite{F-O98}.
\\[5mm]
In approximation problems, the Douglas-Shapiro-Shields
factorization is often used to represent a transfer function.
Since in such problems the lossless factor is only determined up
to a unitary constant matrix, we are interested in a
parametrization of the quotient space of lossless functions by
unitary matrices. To be precise, two lossless functions $G_1(z)$
and $G_2(z)$ in this quotient are equivalent if there exists a
unitary matrix $X$ such that $G_1(z) = X \, G_2(z)$. To deal with
such an equivalence relation, we have the identity
\[
  \cT_{\widehat{\Theta}(w,u, X \, v)}(X \, G) =
  X \, \cT_{\widehat{\Theta}(w,u,v)}(G),
\]
which is satisfied for every lossless $G$ and every unitary $X$.
Therefore, if $G(z)$ has a set of Schur parameters
$(v_1,\ldots,v_n)$ and a unitary matrix $G^{(0)}=D_0$ in a chart
defined by given sequences $(w_1,\ldots,w_n)$ of interpolation
points and $(u_1,\ldots,u_n)$ of direction vectors, then $X \,
G(z)$ has the Schur parameters $(X \, v_1, \ldots, X \, v_n)$ and
the unitary matrix $G^{(0)} = X \, D_0$ in the same chart. An
atlas of the quotient space can therefore be obtained by fixing
the unitary matrix $G^{(0)}$ in each chart instead of letting it
vary. For $w_1,w_2,\ldots,w_n$ in the open unit disk,
$u_1,u_2,\ldots,u_n \in {\mathbb C}^p$ of unit length {\em and}
$D_0$, a fixed $p \times p$ unitary matrix, we define the chart
$(\widetilde{\cV},\widetilde{\varphi})$ by its domain
\[
  \widetilde{\cV} (w_1,w_2,\ldots,w_n,u_1,u_2,\ldots,u_n,D_0) =
  \{ G  \in \cL_n^p \, | \, \|G^{(k)}(1/\overline{w}_k)u_k\|<1,
  G^{(0)}=D_0 \},
\]
and its coordinate map
\[
  \begin{array}{cccl}
  \widetilde{\varphi} : & G & \to & (v_1,v_2,\ldots,v_n)
  \end{array}
\]
where the lossless functions $G^{(k)}$ and the Schur parameter
vectors $v_k$ are again those of Theorem \ref{schuralgo}, in which
$\Theta_k = \widehat{\Theta}(u_k,v_k,w_k)$.
\begin{theorem}
The family $(\widetilde{\cV},\widetilde{\varphi})$ defines a
$C^{\infty}$-atlas on the quotient space $\cL_n^p/\cU_p$.
\end{theorem}
Moreover, a unique balanced realization computed from the Schur
parameters corresponds to every lossless function in the domain
$\widetilde{\cV}$ of a chart, by iterating formula
(\ref{statespacerecursion}) in which $U$ and $V$ are the unitary
matrices given in Theorem \ref{generalcase}. These `Schur balanced
realizations' are really in a (local) canonical form. With each
$(A,B,C,D)$ this method associates (for a given choice of
interpolation points and direction vectors) a unique realization
$(\widetilde{A}, \widetilde{B}, \widetilde{C}, \widetilde{D})$.
Attaching a canonical form with each chart, we obtain a set of
overlapping canonical forms. This generalizes the results of
\cite{HP2000} to the multivariable case and opens up possibilities
for multivariable stable all-pass model reduction and
approximation methods.

Note that, except for particular choices of the interpolation
points and direction vectors, these overlapping canonical forms in
general cannot be characterized in a simple way, and cannot be
computed directly from a given realization. However, if the
interpolation points $w_k$ are the poles of $G(z)$ with the
direction vectors $u_k$ spanning the associated kernels of the
singular matrices $G(1/\overline{w}_k)$, then the Schur vectors
$v_k$ are all zero and the Schur balanced realization is in Schur
form (i.e., the matrix $A$ is triangular). If the interpolation
points $w_k$ are all equal to zero and the vectors $u_k$ are all
equal to $e_1=[1,0,\ldots,0]^*$, then the Schur balanced
realization is positive upper Hessenberg.

The overlapping balanced canonical forms for stable all-pass
systems obtained in the way described above, also give rise to
output-normal canonical forms, resp. input-normal canonical forms,
as follows. Consider an arbitrary (asymptotically) stable system
in minimal state-space form, given by $(A,B,C,D)$. It is not hard
to show that one can find a pair $(\widetilde{B},\widetilde{D})$
such that $(A,\widetilde{B},C,\widetilde{D})$ is stable all-pass.
One can then choose a state-space transformation that brings this
system into one of the balanced canonical forms presented in this
paper, which we denote by $(A_b,B_b,C_b,D_b)$. The resulting
state-space transformation $T$ and the resulting pair
$(A_b,C_b)=(TAT^{-1},CT^{-1})$ are independent of the choice of
$\widetilde{B},\widetilde{D}$.

Therefore applying the same state-space transformation to
$(A,B,C,D)$ one obtains a canonical form $(A_n,B_n,C_n,D_n)$,
where $A_n=A_b$, $C_n=C_b$, $B_n=TB$, $D_n=D$. Obviously the
observability Gramian of this canonical form is the identity
matrix, hence it is an output-normal canonical form. By duality
one can also construct input-normal canonical forms in the same
fashion. These output-normal and input-normal canonical forms are
of importance in model order reduction and system identification
problems. An analogous canonical form in continuous time was used
in \cite{HanzonMaciejowski1996} to construct an algorithm to find
the stable system of order $k$, $k<n$, which is closest in the
$H_2$-norm to a given stable system of order $n$. In system
identification output-normal canonical forms are used in the
context of the so-called separable least squares approach (cf.,
e.g., \cite{BCHV1996}).

The parametrizations described in this paper have been implemented
in a software program named RARL2, dedicated to the rational
approximation of {\em multivariable} discrete-time transfer
functions in $L^2$ norm. Just as the software program Hyperion, it
is based on the algorithm described in \cite{F-O98}. However, in
contrast to Hyperion it has the particularity to deal with
state-space formulations. This state-space approach is of great
interest due to the good numerical behavior of the recursive
construction of balanced realizations presented in this paper.
Presently the main domain of application of RARL2 is the synthesis
and identification of filters from partial frequency data.



{\bf Acknowledgements.}
This research was supported in parts by the Van
   Gogh program VGP 61-431. The authors acknowledge financial support in part by the
 European Commission through the program Training and Mobility of
 Researchers - Research Networks and through project System
 Identification (FMRX CT98 0206) and acknowledge contacts with the
 participants in the European Research Network System
 Identification (ERNSI).

%


\begin{thebibliography}{99}
\bibitem{ABG94}{D.\ Alpay, L.\ Baratchart and A.\ Gombani,
        On the differential structure of matrix-valued rational inner
        functions, {\em Operator Theory: Adv.\ and Appl.}, {\bf 73},
        30--66, 1994.}
\bibitem{BC}{M.\ Bakonyi and T.\ Constantinescu,
        {\it Schur's algorithm and several applications},
        Longman Sci.\ Tech., Harlow, 1992.}
\bibitem{BGR}{J.A.\ Ball, I.\ Gohberg and L.\ Rodman,
        Interpolation of rational matrix functions,
        {\em Operator Theory: Advances and Applications}, {\bf 45},
        Birkh\"{a}user, 1990.}
\bibitem{Bar90}{L.\ Baratchart,
        On the topological structure of inner functions and its use in
        identification, in: {\em Analysis of Controlled Dynamical Systems,
        Lyon, France, 1990, Progress in Systems and Control Theory}, Vol.\
        {\bf 8}, Birkh\"{a}user, 51--59, 1990.}
\bibitem{BOiuf}{{L.\ Baratchart and M.\ Olivi},
        {Inner-unstable factorization of stable rational transfer functions}
        in: {\em Progress in System and Control Theory:
        Modeling, Estimation and Control of Systems with Uncertainty}, Vol.\
        {\bf 10}, {Birkh\"{a}user}, {22--39}, {1991}.}
\bibitem{BCHV1996}{ J.\ Bruls, C.T.\ Chou, B.R.J.\ Haverkamp and M.\ Verhaegen,
        {\em Linear and non-linear system identification using separable
        least-squares}, Technical University Delft, Dept. Electrical Engineering,
        Systems- and Control Laboratory, Report TUD/ET/SCE96.009, 1997.}
\bibitem{DSS70}{R.\ Douglas, H.\ Shapiro and A.\ Shields,
        Cyclic vectors and invariant subspaces for the backward shift
        operator, {\em Annales de l'Institut Fourier,
        (Grenoble)}, {\bf 20}, 37--76, 1970.}
\bibitem{Dym}{H.\ Dym, $J$-contractive matrix functions, reproducing kernel
        spaces and interpolation, {\em {CBMS} lecture notes}, {\bf 71},
        American Mathematical Society, Rhode Island, 1989.}
\bibitem{Dym89}{H.\ Dym, On reproducing kernel spaces, $J$-unitary matrix
        functions, interpolation and displacement rank,
        {\em Operator Theory~:~Advances and Applications},
        Vol.\ {\bf 41}, 173--239, 1989.}
\bibitem{F-O98}{P.\ Fulcheri and M.\ Olivi,
        Matrix Rational $H^2$-Approximation: A Gradient Algorithm
        Based On Schur Analysis,  {\em SIAM Journal on Control and
        Optimisation}, Vol.\ {\bf 36}, No.\ 6, 2103--2127, 1998.}
\bibitem{Garnett}{J.\ Garnett,
        {\em Bounded Analytic Functions}, Academic Press, 1981.}
\bibitem{Getal}{Y.\ Genin, P.\ Van Dooren, T.\ Kailath, J.-M.\ Delosme and
        M.\ Morf, On $\Sigma$-lossless transfer functions and related questions,
        {\em Linear Algebra and its Applications}, {\bf 50}, 251--275, 1983.}
\bibitem{HanzonMaciejowski1996}{B.\ Hanzon and J.M.\ Maciejowski,
        Constructive algebra methods for the $L_2$ problem for stable linear
        systems, {\em Automatica}, Vol.\ {\bf 32}, No.\ 12, 1645--1657, 1996.}
\bibitem{HP2000}{B.\ Hanzon and R.L.M.\ Peeters,
        Balanced Parametrizations of Stable SISO All-Pass Systems in
        Discrete-Time, {\em MCSS}, {\bf 13}, 240--276, 2000.}
\bibitem{H-O97}{B.\ Hanzon and R.J.\ Ober,
        Overlapping block-balanced canonical forms and parametrizations:
        the SISO case,
        {\em SIAM J. Control and Optimization},
        {\bf 35}, 228--242, January 1997.}
\bibitem{H-O98}{B.\ Hanzon and R.J.\ Ober,
        Overlapping block-balanced canonical forms for various classes
        of linear systems, {\em Linear Algebra and its Applications},
        {\bf 281}, 171--225, 1998.}
\bibitem{Kailath}
        {T.\ Kailath, {\em Linear Systems}, Prentice-Hall, 1980.}
\bibitem{Kimura}
        {H.\ Kimura, {\em Chain-Scattering Approach to $H^\infty$-Control},
          Birkh\"auser, 1997.}
\bibitem{L-O98}{J.\ Leblond and M.\ Olivi,
        Weighted $H^2$ Approximation of Transfer Functions, {\em MCSS},
        {\bf 11}, 28--39, 1998.}
\bibitem{Moore}{B.C.\ Moore,
        Principal component analysis in linear systems: Controllability,
        observability, and model reduction,
        {\em IEEE Trans. Autom. Control} {\bf AC-26}, 17--32, 1981.}
\bibitem{Obe87}{R.J.\ Ober,
        Balanced realizations: canonical form, parametrization, model
        reduction, {\em Int.\ J.\ Control}, {\bf 46}, 643--670, 1987.}
\bibitem{OberMunich1987}{R.J.\ Ober,
        Asymptotically stable allpass transfer functions:
        canonical form, parametrization and realization, in:{\em Proceedings of the
        IFAC World Congress,} Munich, 1987.}
\bibitem{Pot60}{V.P.\ Potapov,
        The multiplicative structure of $J$-contractive matrix
        functions, {\em Amer.\ Math.\ Soc.\ Transl.} (2), {\bf 15},
        131--243, 1960.}
\bibitem{Pot88}{V.P.\ Potapov,
        Linear Fractional Transformations of Matrices, {\em Amer.\
        Math.\ Soc.\ Transl.} (2), {\bf 138}, 21--35, 1988.}
\bibitem{PS82}{L.\ Pernebo, L.M.\ Silverman,
        Model reduction via balanced state space representations.
        {\em IEEE Trans. Autom. Control}, {\bf AC-27}, 382--387, 1982.}
\bibitem{Young}{N.J.\ Young,
        Linear fractional transformations in rings and modules,
        {\em Linear Algebra Appl.}, {\bf 56}, 251--290, 1984.}
\end{thebibliography}
\end{document}